\documentclass{article} 
\usepackage[english]{babel}

\usepackage[utf8]{inputenc}
\usepackage{amsmath,amsfonts,amssymb,amscd}
\usepackage{amsthm}
\usepackage{booktabs}
\usepackage{multirow}
\usepackage{alltt}
\usepackage{algorithm}
\usepackage{algorithmic}

\usepackage{color}
\usepackage[pdftex]{graphicx}
\usepackage{epstopdf}

\usepackage{tikz}
\usetikzlibrary{arrows,calc,decorations.markings,intersections,positioning,shapes.arrows}

\usepackage{caption}
\usepackage{subcaption}

\usepackage{cite}

\usepackage{url}
\usepackage{todonotes}

\newcommand{\R}{\mathbb{R}}

\newcommand{\half}{\frac{1}{2}}

\newcommand{\norm}[1]{\left\|#1\right\|}

\newcommand{\DD}{D^{2}}
\newcommand{\Dt}{\Delta t}

\theoremstyle{plain}
\newtheorem{theorem}{Theorem}[section]

\theoremstyle{definition}

\newtheorem{remark}[theorem]{Remark}

\usepackage{pgfplots}
\usepackage{pgfplotstable}
\pgfplotsset{compat=newest}
\usepackage{booktabs}
\usepackage{makecell}

\pgfplotstableset{
every head row/.style={
before row=\toprule,after row=\midrule},
every last row/.style={
after row=\bottomrule},
}

\pgfplotstableset{timingsStyle/.append style={  %
create on use/efficiency/.style = {
   create col/expr={\thisrow{speedup}/\thisrow{nproc}*100},
},
columns={nproc,walltime,speedup,efficiency, final_rel_L2_err_ctrl},
columns/nproc/.style={
column name={\#CPUs},
int trunc,
},
columns/walltime/.style={
column name=\makecell{time [s]},
 fixed,zerofill, precision=1,
 dec sep align,
 postproc cell content/.append style={
 }
},
columns/speedup/.style={
column name={speedup},
 std, precision=1,
 string replace={0}{}, 
 empty cells with={--}
},
columns/efficiency/.style={
column name={efficiency [\%]},
 fixed, zerofill, precision=1,
 string replace={0}{}, 
 empty cells with={--},
},
columns/final_rel_L2_err_ctrl/.style={
column name={rel. err. to $u_\text{exact}$},
 sci, zerofill, precision=1,
},
}
}

\pgfplotstableset{timingsStyleHeat/.append style={  %
create on use/efficiency/.style = {
   create col/expr={\thisrow{speedup}/\thisrow{nproc}*100},
},
columns={nproc, walltime, speedup, efficiency},
columns/nproc/.style={
column name={\#CPUs},
int trunc,
},
columns/walltime/.style={
column name=\makecell{time [s]},
 fixed,zerofill, precision=1,
 dec sep align,
 postproc cell content/.append style={
 }
},
columns/speedup/.style={
column name={speedup},
 std, precision=1,
 string replace={0}{}, 
 empty cells with={--}
},
columns/efficiency/.style={
column name={efficiency [\%]},
 fixed, zerofill, precision=1,
 string replace={0.0}{--}, 
 string replace={0}{--}, 
 empty cells with={--},
},
}
}

\pgfplotstableset{warmstartTimingsStyle/.append style={  %
create on use/efficiency/.style = {
   create col/expr={\thisrow{speedup}/\thisrow{nproc}*100},
},
columns={nproc,walltime,speedup,speedup_to_cold,efficiency, final_rel_L2_err_ctrl},
columns/nproc/.style={
column name={\#CPUs},
int trunc,
},
columns/walltime/.style={
column name=\makecell{time [s]},
 fixed,zerofill, precision=1,
 dec sep align,
 postproc cell content/.append style={
 }
},
columns/speedup/.style={
column name={total speedup},
 std, precision=1,
 string replace={0}{}, 
 empty cells with={--}
},
columns/speedup_to_cold/.style={
column name=\makecell{speedup\\to cold},
 std, precision=1,
 string replace={0}{}, 
 empty cells with={--}
},
columns/efficiency/.style={
column name={efficiency [\%]},
 fixed, zerofill, precision=1,
 string replace={0}{}, 
 empty cells with={--},
},
columns/final_rel_L2_err_ctrl/.style={
column name={rel. err. to $u_\text{exact}$},
 sci, zerofill, precision=1,
},
}
}

\pgfplotstableset{warmstartTimingsStyleHeat/.append style={  %
create on use/efficiency/.style = {
   create col/expr={\thisrow{speedup}/\thisrow{nproc}*100},
},
columns={nproc,walltime,speedup,speedup_to_cold,efficiency},
columns/nproc/.style={
column name={\#CPUs},
int trunc,
},
columns/walltime/.style={
column name=\makecell{time [s]},
 fixed,zerofill, precision=1,
 dec sep align,
 postproc cell content/.append style={
 }
},
columns/speedup/.style={
column name={total speedup},
 std, precision=1,
 string replace={0}{}, 
 empty cells with={--}
},
columns/speedup_to_cold/.style={
column name=\makecell{speedup\\to cold},
 std, precision=1,
 string replace={0}{}, 
 empty cells with={--}
},
columns/efficiency/.style={
column name={efficiency [\%]},
 fixed, zerofill, precision=1,
 string replace={0.0}{}, 
 string replace={0}{}, 
 empty cells with={--},
  postproc cell content/.append code={%
    \ifnum\pgfplotstablerow=0
    \pgfkeyssetvalue{/pgfplots/table/@cell content}{--}%
    \fi
     },
},
}
}

\pgfplotsset{select coords between index/.style 2 args={
    x filter/.code={
        \ifnum\coordindex<#1\fi
        \ifnum\coordindex>#2\fi
    }
}}

\title{An Efficient Parallel-in-Time Method for Optimization with Parabolic PDEs\thanks{Submitted to the editors January 17, 2019}}

\author{Sebastian Götschel\thanks{Zuse Institute Berlin, Berlin, Germany, {goetschel@zib.de}} \and Michael L. Minion\thanks{Lawrence Berkeley National Laboratory, Berkeley, CA, USA, {mlminion@lbl.gov}}}

\date{July 10, 2019}

\begin{document}

\maketitle

\begin{abstract}
To solve optimization problems with parabolic PDE constraints, 
often methods working on the reduced
objective functional are used. They are
computationally expensive due to the necessity of
solving both the state equation and a
backward-in-time adjoint equation to evaluate the reduced gradient
in each iteration of the
optimization method. In this study, we investigate the use of
the parallel-in-time method PFASST in the setting of PDE-constrained optimization.
In order to develop an efficient fully time-parallel algorithm
we discuss different options for applying PFASST to adjoint gradient computation, 
including the possibility of doing PFASST
iterations on both the state and adjoint equations simultaneously.
We also explore the additional gains in efficiency from reusing information from
previous optimization iterations when solving each equation.
Numerical results for both a linear and a non-linear
reaction-diffusion optimal control problem
demonstrate the parallel speedup and efficiency of different approaches.
\end{abstract}




\section{Introduction} 
\label{section:introduction}

Large-scale PDE-constrained optimization problems occur in a multitude of applications, for example in solving inverse problems for non-destructive testing of materials and structures~\cite{GoetschelEtAl2016} or in individualized medicine~\cite{FinsbergEtAl2018}. More recently, the training of certain deep neural networks in machine learning, e.g, for image recognition or natural language processing, has been formulated as a dynamic optimal control problem~\cite{HaberRuthotto2017,RuthottoHaber2018}.
Algorithms for the solution of PDE-constrained optimization problems are computationally extremely demanding as they require to numerically solve multiple partial differential equations during the iterative optimization process. This is especially challenging for transient problems, where the solution of the associated optimality system requires information about the discretized variables on the whole space-time domain. For the solution of such optimization problems, methods working on the reduced objective functional are often employed to avoid a full spatio-temporal discretization of the problem. The evaluation of the reduced gradient then requires one solve of the state equation forward in time, and one backward-in-time solve of the adjoint equation. In order to tackle realistic applications, it is not only essential to devise efficient discretization schemes for optimization, but also to use advanced techniques to exploit computer architectures and decrease the time-to-solution, which otherwise is prohibitively long.
%
One approach is to utilize the number of CPU cores of current and future many-core high performance computing systems by parallelizing the PDE solution method. In addition to well-established methods for parallelization in the spatial degrees of freedom, parallel-in-time methods have seen a growing interest in the last 15 years. Research into time parallelism dates back at least to the 1960s, 
and in 2001 the introduction of parareal by Lions et al.~\cite{LionsMadayTurinici2001} sparked new research into time-parallel methods. 
The methods in this paper are based on the parallel full approximation scheme in space and time (PFASST) introduced by Emmett and Minion~\cite{EmmettMinion2012}.
We will not attempt a thorough review of the field here, and the interested reader is  
encouraged to consult the survey article~\cite{Gander2015} for an overview of competing approaches. 

More recently, the application of space-time parallel
methods to the solution of optimization problems governed by PDEs
has become an active research area, with approaches including multiple shooting (e.g.,~\cite{Heinkenschloss2005}
and the references therein), Schwarz methods~\cite{BarkerStoll2015,GanderKwok2016}, and the application of
parareal preconditioners~\cite{MathewEtAl2010,Ulbrich2015}. A time-parallel
gradient type method is presented in~\cite{DengHeinkenschloss2016}. There 
the time interval of interest is subdivided into time steps,
which are solved in parallel using quantities from the previous
optimization iteration as input. This leads to jumps in the solutions
of state and adjoint equation such that these equations are not
satisfied during optimization. While they report excellent speedups and linear scaling up to 50 processors and show convergence if sufficiently small step sizes for 
updating the control are used, it is unclear how to
automatically select such a step size. Alternatively, space-time parallel
multigrid methods are applied to adjoint gradient computation 
and simultaneous optimization~\cite{GuentherGaugerSchroder2017,GuentherGaugerSchroder2018} within the XBraid software library \cite{xbraid-package}. 
XBraid provides a non-intrusive framework adding time-parallelism to existing serial time stepping codes, 
and using simultaneous instead of reduced space optimization, a speedup of 19 using 256 time processors has been reported. The same method is also applied to perform layer-parallel training of neural networks~\cite{GueRuSchroeCyrGau2018}.

In this paper we employ PFASST to provide a fully time-parallel reduced-space gradient- or nonlinear conjugate gradient method that allows using the usual line search criteria, e.g., the strong Wolfe conditions, for step size selection to guarantee convergence, and is thus non-intrusive with respect to the optimization algorithm. On the other hand this is an intrusive approach concerning the PDE solvers as it requires using multilevel spectral deferred correction methods as time steppers. The implementation effort is mitigated by the availability of libraries like LibPFASST~\cite{libpfasst} or dune-PFASST~\cite{dune-pfasst}. 
While the basic ideas of the approach are outlined in the short paper~\cite{GoetschelMinion2018}, here we provide more details, and develop additional approaches to increase speedup and efficiency.
We demonstrate that using PFASST brings additional benefits, like flexibility in treating nonlinearities, and, more importantly, the possibility of warm starting the PDE-solutions required during the optimization.
The remainder of the paper is organized as follows. In Sect.~\ref{section:oc} we present the optimization problem and review optimality conditions as well as adjoint gradient computation. The PFASST method is introduced in Sect.~\ref{section:PFASST}; it is used to derive parallel-in-time methods for solving optimization problems with parabolic PDEs in Sect.~\ref{section:PFASSTOC}. Finally, in Sect.~\ref{section:results} we present numerical examples, followed by a discussion of results and future improvements in Sect.~\ref{section:outlook}.

\section{Adjoint gradient computation for optimization with parabolic PDEs}
\label{section:oc}

Here we briefly summarize the mathematical approach to parabolic PDE-constrained optimization problems. For more details and generalizations, we refer to, e.g.,~\cite{HinzeEtAl2009}.

We consider optimization problems of the form
\begin{equation}\label{eq:ocp}
 \min_{y\in Y, u\in U} J(y,u)\ \text{subject to}\ c(y,u) = 0,
\end{equation}
with the equality constraint $c : Y\times U \rightarrow Z^\star$ being a parabolic PDE on Hilbert spaces $Y, U, Z$. $Z^\star, U^\star$ denote the dual spaces of $Z$ and $U$, respectively; the dual pairing between a space $X$ and its dual is denoted by $\langle\cdot,\cdot\rangle_{X^\star,X}$, and $(\cdot,\cdot)_X$ denotes the scalar product on $X$. We drop the subscripts like $\cdot_X$ if the involved spaces are clear from the context.
In the present setting, the constraint
$c(y,u) = 0$ means that the \emph{state} $y$ satisfies a PDE where the \emph{control} $u$ is a specific forcing term, occurring, e.g., as a source term, in the boundary conditions, or as some other parameter.

To derive optimality conditions, we assume that there exists a unique solution $y = y(u) \in Y$ of the state equation $c(y,u) = 0$ for each control $u \in U$. We additionally assume that $c_y(y,u) : Y\rightarrow Z^\star$ is continuously invertible. Then, by the implicit function theorem (see, e.g.,~\cite[Section 4.7]{Zeidler1986}), the control-to-state mapping is continuously differentiable, and the derivative $y'(u)$ is given by the solution of
\begin{equation}
 c_y(y,u) y'(u) + c_u(y,u) = 0.
\end{equation}
Subscripts like $c_u()$ denote the partial derivatives with respect to the indicated variable.
By inserting $y(u)$ into the optimization problem~\eqref{eq:ocp} we arrive at the reduced problem
\begin{equation}\label{eq:reduced_ocp}
 \min_{u\in U} j(u) := J(y(u),u).
\end{equation}
In this unconstrained setting, the following simple first-order necessary optimality condition holds. If $u^\star \in U$ is a local solution of the reduced problem~\eqref{eq:reduced_ocp} it is a zero of the reduced derivative, $j'(u^\star) = 0.$ If the reduced functional $j$ is convex, this condition is also sufficient.
If we allow control constraints, i.e.,~demand $u\in U_\text{ad}$ with $U_\text{ad} \subset U$ non-empty, convex and closed, the optimality condition changes to the variational inequality for the local minimizer $u^\star \in U_\text{ad}$
\begin{equation}
 \langle j'(u^\star),u-u^\star\rangle_{U^\star,U} \geq 0 \quad \forall u \in U_\text{ad}.
\end{equation}

To formally derive a representation for the reduced gradient, we define the Lagrange functional $\mathcal{L} : Y \times U \times Z \rightarrow \R$,
\begin{equation}\label{eq:lagrange}
 \mathcal{L}(y,u,p) = J(y,u) + \langle p, c(y,u)\rangle_{Z,Z^\star},
\end{equation}
where in the present context, the Lagrange multiplier $p \in Z$ is referred to as the \emph{adjoint}. 
Clearly, inserting $y = y(u)$ into~\eqref{eq:lagrange}, we get $j(u) = \mathcal{L}(y(u),u,p)$ for arbitrary $p \in Z$.
Differentiation in direction $\delta u\in U$ yields
\begin{equation}\label{eq:directionalderivative}
 \langle j'(u), \delta u\rangle_{U^\star,U} = \langle \mathcal{L}_y(y(u),u,p), y'(u)\delta u\rangle_{Y^\star,Y} + \langle \mathcal{L}_u(y(u),u,p), \delta u\rangle_{U^\star,U}.
\end{equation}
Choosing $p = p(u)$ such that the adjoint equation
\begin{equation}\label{eq:adjoint}
  c_y(y(u),u)^\star p(u) = -  J_y(y(u),u)
\end{equation}
is fulfilled gives
\begin{equation}\label{eq:Ly0}
 \mathcal{L}_y(y(u),u,p(u)) = J_y(y(u),u) + c_y(y(u),u)^\star p(u) = 0.
\end{equation}
Inserting this into~\eqref{eq:directionalderivative}, the first term on the right hand side vanishes, and we get the reduced derivative $j'(u) \in U^\star$ as
\begin{equation*}
\langle j'(u), \delta u\rangle_{U^\star,U} = \langle \mathcal{L}_u(y(u),u,p), \delta u\rangle_{U^\star,U},
\end{equation*}
i.e.,
\begin{equation}\label{eq:redgrad}
j'(u) = \mathcal{L}_u(y(u),u,p(u)) = J_u(y(u),u) + c_u(y(u),u)^\star p(u).
\end{equation}
In the Hilbert space setting used here, for a given $u\in U$ the reduced gradient $\nabla j(u) \in U$ is then given as the Riesz representative of the reduced derivative $j'(u) \in U^\star$,
i.e.,~via, \[ \bigl( \delta u, \nabla j(u)\bigr)_U = j'(u)\delta u\ \ \forall \delta u\in U.\] 

To be more concrete, we consider a distributed tracking-type objective functional $J(y,u)$ with an additional term penalizing the control cost for $\lambda \ge 0$,
\begin{equation}\label{eq:objective}
J(y,u) = \underbrace{\half \int_0^T \norm{y-y_d}^2_{L^2(\Omega)} \ dt}_{=J^\Omega(y,u)}+ \underbrace{\frac{\lambda}{2} \int_0^T \norm{u}_{L^2(\Omega)} \ dt}_{=J^u(y,u)}
\end{equation}
with $y_d$ representing a desired solution.
The state  $y$ is subject to a  linear or semi-linear parabolic PDE with distributed control
\begin{alignat}{3}
y_t - \kappa \Delta y + f(y) - u &= 0  &\quad&\text{in}\ \Omega\times [0,T] \label{eq:state1} \\
y(0)-y_0 &= 0  &&\text{in}\ \Omega \label{eq:state2}
\end{alignat}
and suitable boundary conditions. For simplicity, we consider a constant scalar diffusivity $\kappa$; $\Delta y$ denotes the Laplacian of $y$. Thus, for sufficiently smooth $f$, we choose $Y=U=Z=L^2(0,T; H)$ with $H=L^2(\Omega)$ in~\eqref{eq:ocp}. Here, $\Omega \subset \R^n$, and $T>0$ denotes a given final time.
The adjoint equation~\eqref{eq:adjoint} for a given state solution $\bar y$ becomes
\begin{alignat}{3}
-p_t - \kappa \Delta p + f_y(\bar y)p  &= -J_y^\Omega(\bar y, u) = -(\bar y-y_d) &\quad&\text{in}\ \Omega\times [0,T] \label{eq:adjoint1} \\
p(T) &= 0 &&\text{in}\ \Omega,\label{eq:adjoint2}
\end{alignat}
with homogeneous boundary conditions of the same type as in the state equation. Generalizations to controls or observations on the boundary are straightforward. Having a terminal cost contribution in the objective, e.g., 
\begin{equation*}
J^T(y,u) = \frac{\sigma}{2} \norm{y(T)-y_d^T}^2_{L^2(\Omega)}
\end{equation*}
leads to the modified terminal condition
\begin{equation}\label{eq:adjoint3}
p(T) = -J^T_y(\bar y,u) = -(y(T)-y_d^T).
\end{equation}

\begin{figure}[ht]
\begin{center}
\begin{tikzpicture}
  \draw[->,ultra thick] (-2,0)--(4,0) node [above] {$c(y,u) = 0$};
  \draw[<-,ultra thick] (-2,-1) node [below] {$c(y,u)^\star p = -J_y(y,u)$} --(4,-1) ;
  \foreach \x in {-1.5,-1,-0.5,0,0.5,1.5,2,2.5,3,3.5}
    \draw[->,red,thick] (\x,-0.1)--(\x,-0.9);
  \draw[red] (1,-0.5)  node {$y(u)$};
  \draw[red,thick] (1,-0.1)--(1,-0.25);
  \draw[->,red,thick] (1,-0.75)--(1,-0.9);
\end{tikzpicture}
\end{center}
\caption{Adjoint gradient computation. The full solution of the forward state equation is required to solve the adjoint equation backward in time.}
\label{fig:adjgradcomp}
\end{figure}
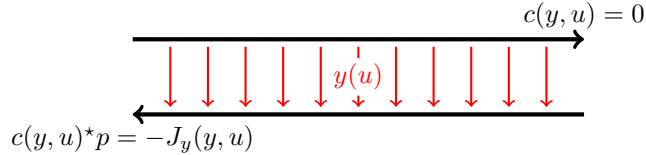

As the adjoint equation~\eqref{eq:adjoint} is backward in time, due to the occurrence of $-J_y(y(u),u)$ as a source term, and---in the nonlinear case---the dependence of $c_y(y(u),u)$ on the state solution $y(u)$, adjoint gradient computation consists of three steps (see also Fig.~\ref{fig:adjgradcomp}):
\begin{enumerate}
 \item Solve the PDE $c(y,u) = 0$ for a give  control $u$,  and store the solution trajectory  $y\in Y$.
 \item Solve the adjoint PDE $c_y(y,u)^\star p = -J_y(y,u)$ for $p\in Z$.
 \item Set $j'(u) = J_u(y,u) + c_u(y,u)^\star p$ and compute the Riesz representation $\nabla j(u)$.
\end{enumerate}

Thus, the computation of the reduced gradient requires the solution of two parabolic PDEs.  
For solving the optimization problem, in this work we consider nonlinear conjugate gradient (ncg)- and steepest descent (sd) methods, as they require only gradient information. 
Hence gradient-based methods for the optimization problems that
iterate over these three steps require many forward and backward PDE solves.  The goal of this paper is to exploit properties of the the PFASST parallel-in-time algorithm to 
reduce the computational complexity of the optimization procedure.

\begin{remark}
 State and adjoint equations can be solved as a coupled system, as is done, e.g., in~\cite{GuettelPearson2016}. In view of parallelizing in time this becomes more difficult, see the brief discussion in Sect.~\ref{subsection:PinTgradient} and Fig.~\ref{fig:combined}; our main focus thus is on gradient computation and using algorithms like steepest descent or nonlinear conjugate gradients to iteratively solve the optimization problem. This can easily be applied to, e.g., control constrained problems, and be extended to second order methods like semi-smooth Newton, or Newton-CG.
\end{remark}

For the gradient-based methods considered here, the optimization iteration proceeds as 
\begin{align}
 u_{k+1} &= u_k + \alpha_k d_k \label{eq:update_u} \\
 d_{k+1} &= -\nabla j(u_{k+1})+ \beta_k d_k, \label{eq:update_d}
\end{align}
where $d_0 = \nabla j(u_0)$, and the choice of $\beta_k$ defines the actual method. Abbreviating $\nabla j(u_k)$ by $g_{k} $, specific variants of ncg include, e.g., 
\begin{alignat*}{3}
\beta_k^\text{FR} &= \frac{(g_{k+1}, g_{k+1} )}{( g_{k}, g_{k})} &\quad&\text{Fletcher-Reeves~\cite{FletcherReeves1964}}, \\
\beta_k^\text{PRP} &= \frac{( g_{k+1},  g_{k+1} -g_{k} ) }{( g_{k}, g_{k} )} &&\text{Polak-Ribiere-Polyak~\cite{PolakRibiere1969,Polyak1969}}, \\
\beta_k^\text{DY} &= \frac{( g_{k+1},g_{k+1} ) }{( d_k,  g_{k+1}- g_{k} )}&&\text{Dai-Yuan~\cite{DaiYuan1999}}.
\end{alignat*}
Using $\beta_k = 0$ yields the usual steepest descent method.

To guarantee convergence under the usual assumptions, the step size $\alpha_k$ is chosen to satisfy the strong Wolfe conditions for ncg,
\begin{align}
j(u_k + \alpha_k d_k) &\leq j(u_k) + c_1 \alpha_k \langle g_k, d_k \rangle \label{eq:armijo}\\
| \langle j'(u_k + \alpha_k d_k), d_k \rangle | & \leq c_2 |(g_k, d_k) | \label{eq:curvature},
\end{align}
$0<c_1<c_2<1$, and just the Armijo condition~\eqref{eq:armijo} for sd~\cite{NocedalWright2006}.

\begin{remark}
In this work we follow the \emph{first optimize, then discretize} approach. To implement the optimization methods based on the optimality conditions of the previous section, we need to discretize the arising parabolic PDEs in time and space. To facilitate using PFASST for time parallelism, this is done using the method of lines approach, so we discretize space first. In the examples in Sect.~\ref{section:results} we use a pseudo-spectral method, but other techniques like finite element, finite difference, or finite volume methods are possible as well. The resulting system of ODEs is then solved using PFASST.
\end{remark}
\section{SDC,  MLSDC, and PFASST}\label{section:PFASST}
\def\qtil {{\tilde{w}}}
\def\qtilE {{\tilde{w_E}}}
\def\qtilI {{\tilde{w_I}}}
\def\tn {{t_n}}
\def\tnp {{t_{n+1}}}
\def\np {{{n+1}}}
\def\kp {{{[k+1]}}}
\def\k {{{[k]}}}
\def\y {{y}}

In this section we give an overview of the parallel-in-time strategy used to solve the state and adjoint equations.
The strategy is mainly based on the PFASST algorithm~\cite{EmmettMinion2012} with some modifications specific to PDE-constrained optimization problems.  PFASST can be thought of
as a time parallel variant of the multi-level spectral deferred correction method (MLSDC) \cite{Speck2013}, which in turn is constructed from the method of 
spectral deferred corrections (SDC) \cite{DuttGreengardRokhlin:2000}.  Since all three of these methods are well established, we give only a brief overview here skewed towards the numerical methods tested in Sect.~\ref{section:results}.

\subsection{SDC Methods} \label{sect:sdc}
Consider the generic ODE over the time  interval $[t_j,t_{j+1}]$ representing one time step
\begin{equation} \label{eq:ODE}
  y'(t)=F(t,y),
\end{equation}
with initial condition  $y(t_j)=y_j$.
Divide the interval $[t_j,t_{j+1}]$ into
$M$  smaller sub-intervals by choosing points $t_m, ~ m=0, \dots, M$ 
corresponding to the Gauss-Lobatto quadrature nodes.

The exact solution of the ODE at each point $t_m$ is given by
\begin{equation} \label{eq:picard}
   y(t_m)  =  y_j + \int_{t_j}^{t_m} F(\tau,y(\tau)) d\tau. 
\end{equation}
The collocation method is defined by the solution of the system of equations derived by applying quadrature rules to \eqref{eq:picard},
\begin{equation} \label{eq:colloc}
   y_m  \approx  y_j + \Delta t \sum_{i=0}^{M} w_{m,i} F(t_i,y_i), 
\end{equation}
where $\Delta t = t_{j+1}-t_j$ and $w_{m,i}$ are the quadrature weights,
\begin{equation*}
w_{m,i} = \frac{1}{\Delta t} \int_{t_j}^{t_m} \ell_i(s) ds,\ m=0,\dots,M,\ i=0,\dots,M,
\end{equation*}
with Lagrange polynomials $\ell_i$ defined by the quadrature nodes $(t_m)_{m=0,\dots,M}$ (here Lobatto IIIA).
Equation~\eqref{eq:colloc} is a (typically nonlinear) system that is $M$ times larger than that of
a single-step implicit method like backward Euler.  It is also equivalent to a fully implicit Runge-Kutta method with the values  $w_{m,i}$ corresponding to the matrix in the Butcher tableaux. Such methods have good stability properties and have formal order of accuracy $2M$ for $M+1$ Lobatto quadrature nodes. (See e.g. \cite{Hairer1991} for a more detailed discussion of collocation and implicit Runge-Kutta methods).

Spectral deferred correction methods (SDC) can be considered as a fixed point iteration to solve the collocation formulation~\eqref{eq:colloc}.  Each  SDC iteration or {\it sweep} consists of stepping through the nodes $t_m$ and updating the solution at that node $y_m^\k$. A typical version of SDC using implicit sub-stepping takes  the form (with $\k$ denoting iteration)
  \begin{equation} \label{eq:sdcsweep}
      y^{\kp}_m = \y_j+ \Delta t \sum_{i=0}^{m} \qtil_{m,i} F(t_i,y^\kp_i) + \Delta t \sum_{i=0}^{M} (w_{m,i}-\qtil_{m,i}) F(t_i,y^\k_i). 
  \end{equation}
Here the values $\qtil_{m,i}=0$ for $i>m$, hence each of the substeps has the same computational complexity of a single backward Euler step.  

One attractive aspect of SDC methods is that they can easily be extended to cases in which \eqref{eq:ODE} can be split into stiff and nonstiff components.  Semi-implicit or Implicit-Explicit (IMEX) methods that split the equation into stiff and nonstiff terms
first appeared in \cite{Minion2003a}.  So called {\it multi-implicit} or MISDC methods that treat two implicit terms in an operator splitting approach are introduced in 
\cite{Bourlioux2003}.  Such splittings are explored in the numerical tests in Sect.~\ref{subsection:nagumo} by using MISDC variants to treat the nonlinearity of the state equation.
In the numerical examples, the values of $\qtil_{m,i}$ for implicit terms are chosen following the $LU$ factorization method of Weiser \cite{Weiser2013}, while those for explicit terms correspond to the usual forward Euler substepping.

\subsection{MLSDC Methods} \label{sect:mlsdc}

Higher-order SDC methods require a relatively large number of function evaluations (explicit and/or implicit) per time step.  One method to reduce the cost of these iterations is to employ a multi-level formulation of iterations where SDC sweeps are done on a hierarchy of discretization levels. In \cite{Speck2013}, so called multi-level SDC (MLSDC) methods are studied where the levels are differentiated by the spatial and/or temporal order and resolution as well as the tolerance of implicit solvers. SDC sweeps are scheduled like V-cycles in multigrid, and coarse level problems are modified with a term analogous to a space-time full approximation scheme (FAS) correction term. In the current study, only the number of spatial degrees of freedom and the number of SDC quadrature nodes is varied on MLSDC or PFASST levels.
MLSDC is the basic building block for the parallel-in-time method PFASST discussed next.

\subsection{PFASST} \label{sect:pfasst}

The parallel full approximation scheme in space and time (PFA\-SST)~\cite{EmmettMinion2012} is, as the name suggests, a method for exploiting both spatial and  temporal parallelism in a manner similar to FAS multigrid methods for nonlinear problems (see \cite{Bolten2017} for a multigrid perspective of PFASST).  As mentioned above, PFASST can be considered to be a pipelined version of the MLSDC method, with each time step being assigned to a separate processor (or groups of processors), such that MLSDC sweeps are performed on multiple time steps in parallel. ``Pipelined'' here refers to the idea that each processor begins a step of the algorithm as soon as new initial conditions are passed forward in time from the previous processor~\cite{Minion:2010}.

\begin{figure}
	\begin{center}	
		\includegraphics[width=0.8\textwidth]{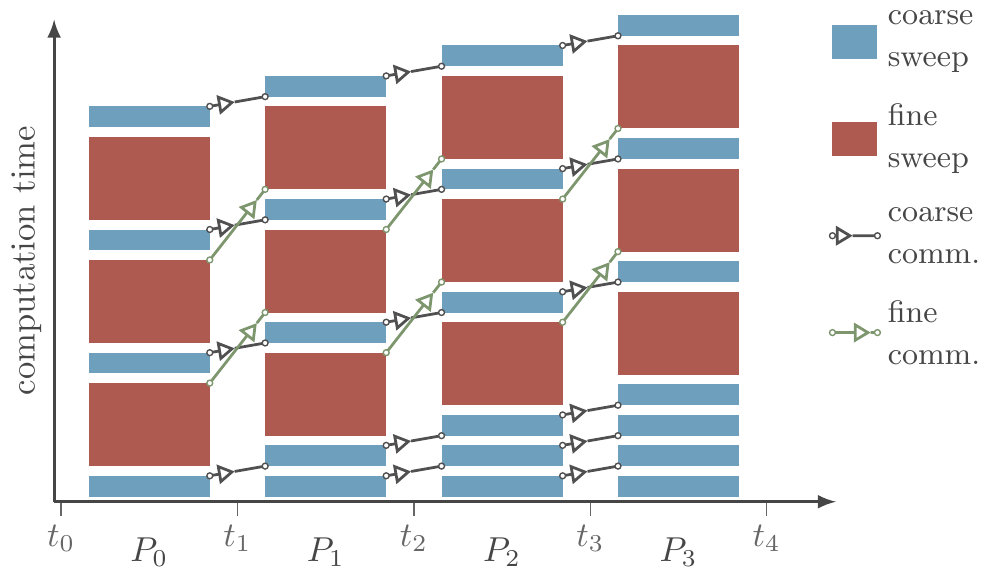}
    \end{center}
    \caption{Generic two-level PFASST scheduling for the forward solution of an ODE~\eqref{eq:ODE}. The time steps  $[t_j,t_{j+1}]$ are distributed among processors $P_j$, which perform coarse SDC sweeps sequentially, and sweeps on the fine level in parallel. The coarse sweeps at the bottom illustrate the predictor, which for PDE-constrained optimization requires communication. Picture created using pfasst-tikz (\protect\url{https://github.com/f-koehler/pfasst-tikz}).}
    \label{fig:normalpfasst}
\end{figure}

As all processors handling time steps $[t_j,t_{j+1}]$ need an initial condition $y_j$, which for $j>0$ is not known, PFASST starts with a \emph{predictor} phase.
Usually, this is done by integrating the equations at the coarsest level in serial (with or without the FAS correction term) with a low-order method. Alternatively, this step can be conducted without coarse level communication since after distributing the initial value $y_0$ to all processors, the $j$th processor can compute $j$ time steps sequentially rather than wait for the
first $j-1$ time steps to be completed before beginning (\emph{burn-in}). 
Note, however, in the context of PDE-constrained optimization, the latter is not possible, as each processor only has data (control, desired state) for the time steps it actually computes, and not for the other time intervals.
Instead, except for the first iteration of optimization method, one can initialize all processors using the solution from the previous optimization iteration.  We refer to this procedure as {\it warm starting} PFASST iterations. In this case, each processor will start doing SDC sweeps with communication at the coarse level propagating the solution forward in time to receive and updated initial value. The predictor step is represented at the bottom of Fig. \ref{fig:normalpfasst}, with warm starting discussed further in Sect. \ref{sect:warm}.

Having finished the predictor phase, each processor performs MLSDC sweeps, where after each SDC sweep on each level of the MLSDC hierarchy, the solution is communicated forward in time to the next processor.  This communication overlaps with computation except at the coarsest level (see Fig.~\ref{fig:normalpfasst} and~\cite{Emmett2012}). A more detailed description can be found in~\cite{EmmettMinion2012}.

\section{Time-parallel PDE-constrained Optimization}
\label{section:PFASSTOC}

In this section, we discuss three components to produce an efficient fully time-parallel method for PDE-constrained optimization. First we briefly describe parallelizing the outer optimization loop, before coming to the more involved parallel-in-time computation of the reduced gradient. Finally, as the time integrators for solving state and adjoint equations are iterative, we discuss using previous solutions to \emph{warm start} the time integration.

For parallelization in time, the time domain $[0,T]$ is subdivided into $N$ time steps $0 = t_0 < \dots < t_N = T$, which are distributed on $R$ processors.
For examples with relatively fewer time steps, the shortest computation times may result from choosing $R=N$.  On the other hand, when the number of processors
used is smaller than the total number of time steps, (e.g., in the strong scaling studies included below), the PFASST algorithm is applied sequentially to blocks of
time steps where $N$ is an integer multiple of $R$.
Note that here we consider only temporal parallelism. Usually, this is used in combination with spatial parallelism, and \emph{multiplies} the speedup gained
from spatial parallelism. This means that the $R$ processors used for parallelizing in time can be considered as groups of processors with additional spatial parallelism. The best strategy for distributing processors between time and space parallelism is in general problem- and machine-dependent and not considered here.

\subsection{Parallelization of the optimization loop}

Time-parallelization of the optimization loop~\eqref{eq:update_u},~\eqref{eq:update_d} requires parallel-in-time computation of the reduced gradient, as well as evaluation of inner products.
In the present optimal control setting for the prototype problem~\eqref{eq:objective}--\eqref{eq:adjoint2} with the usual $L^2$-regularity in time, $( \cdot, \cdot )$ denotes the inner product in $L^2(0,T; H)$, so
\begin{equation}
(v, w) = \int_0^T (v(t), w(t))_H \ dt.
\end{equation}
The evaluation of these inner products, and thus computation of $\beta_k$ as well as evaluation of the objective function $j(u_k), j(u_k+\alpha_k d_k)$ can easily be done time-parallel, as
\begin{equation}
(v, w) = \sum_{i=0}^N  \int_{t_i}^{t_{i+1}} (v(t), w(t))_H.
\end{equation}
After discretization, each processor evaluates the discrete spatial inner product and integrates in time only for its own time steps. In the numerical examples below, evaluation of the time integral is done using a simple trapezoidal rule; of course using other quadrature rules, like re-using the spectral quadrature matrices provided by PFASST, is possible as well. In terms of communication, for each inner product only one scalar value per processor has to be transmitted to one master processor collecting the results. The essential ingredient for a fully parallel optimization method is the computation of the reduced gradient $\nabla j(u_k)$, which is discussed next.

\subsection{Time-parallel adjoint gradient computation}
\label{subsection:PinTgradient}

Here we discuss three different ways in which PFASST can be used to enable a time-parallel computation of the reduced gradient $\nabla j(u_k)$.
In the simplest case, PFASST is used to first solve the state equation, then to solve the adjoint equation afterwards (Sect.~\ref{subsection:vanilla}). Alternatively, if $R=N$, the adjoint can be solved simultaneously with the state equation (Sect.~\ref{subsection:simul}). 
A third variant is inspired by the paraexp method~\cite{GanderGuettel2013}. In this variant, we make use of the linearity of the adjoint equation to split the adjoint solve into an inhomogeneous equation with homogeneous terminal conditions on each time step (without communication), and a subsequent propagation of the correct terminal conditions backward-in-time (Sect.~\ref{subsection:mixed}).

\subsubsection{First state, then adjoint}\label{subsection:vanilla}

\begin{figure}
\begin{center}
\begin{tikzpicture}
\node[shape=circle,draw,inner sep=2pt] at (-1,0) (char){1};
\draw[thick, black] (0,0) -- (4.25,0);
\draw[thick, black, ->] (4.75,0) -- (9,0) node[right] {$t$};
\draw (0.5,-3pt) -- (0.5,3pt) node[below = 1mm] {$t_{j}$};
\draw (1.375,-1.5pt) -- (1.375,1.5pt) ;
\draw (2.25,-1.5pt) -- (2.25,1.5pt) ;
\draw (3.125,-1.5pt) -- (3.125,1.5pt) ;
\draw (4.0,-3pt) -- (4.0,3pt) node[below = 1mm] {$t_{j+1}$};
\draw (5.0,-3pt) -- (5.0,3pt) node[below = 1mm] {$t_{j+1}$};
\draw[thick, black, ->] (0.65, 5pt) -- (3.85, 5pt);
\draw[thick, red, ->] (4.0, 5pt) to [out=75, in=105] (5.0, 5pt); 
\draw[thick, black, ->] (5.15, 5pt) -- (8.35, 5pt);
\draw (5.875,-1.5pt) -- (5.875,1.5pt) ;
\draw (6.75,-1.5pt) -- (6.75,1.5pt) ;
\draw (7.625,-1.5pt) -- (7.625,1.5pt) ;
\draw (8.5,-3pt) -- (8.5,3pt) node[below = 1mm] {$t_{j+2}$};
\node[shape=circle,draw,inner sep=2pt] at (-1,-0.75) (char){2};
\draw[thick, black] (0,-0.75) -- (4.25,-0.75);
\draw[thick, black, ->] (4.75,-0.75) -- (9,-0.75) node[right] {$t$};
\draw ($(0.5,-0.75)+(0,-3pt)$) -- ($(0.5,-0.75)+(0,+3pt)$) node[below = 1mm] {$t_{j}$};
\draw ($(1.375,-0.75)+(0,-1.5pt)$) -- ($(1.375,-0.75)+(0,1.5pt)$);
\draw ($(2.25,-0.75)+(0,-1.5pt)$) -- ($(2.25,-0.75)+(0,1.5pt)$) ;
\draw ($(3.125,-0.75)+(0,-1.5pt)$) -- ($(3.125,-0.75)+(0,1.5pt)$);
\draw ($(4.0,-0.75)+(0,-3pt)$) -- ($(4.0,-0.75)+(0,3pt)$) node[below = 1mm] {$t_{j+1}$};
\draw ($(5.0,-0.75)+(0,-3pt)$) -- ($(5.0,-0.75)+(0,3pt)$) node[below = 1mm] {$t_{j+1}$};
\draw[thick, black, <-] ($(0.65, -0.75)+(0,-5pt)$) -- ($(3.55,-0.75)+(0,-5pt)$);
\draw[thick, red, <-] ($(4.0, -0.75)+(0,-10pt)$) to [out=-75, in=-105] ($(5.0, -0.75)+(0,-10pt)$); 
\draw[thick, black, <-] ($(5.25, -0.75)+(0,-5pt)$) -- ($(8.15,-0.75)+(0,-5pt)$);
\draw ($(5.875,-0.75)+(0,-1.5pt)$) -- ($(5.875,-0.75)+(0,1.5pt)$) ;
\draw  ($(6.75,-0.75)+(0,-1.5pt)$) -- ($(6.75,-0.75)+(0,1.5pt)$) ;
\draw  ($(7.625,-0.75)+(0,-1.5pt)$) -- ($(7.625,-0.75)+(0,1.5pt)$) ;
\draw ($(8.5,-0.75)+(0,-3pt)$) -- ($(8.5,-0.75)+(0,3pt)$) node[below = 1mm] {$t_{j+2}$};
\end{tikzpicture}
\end{center}
\caption{Two time steps of the first state, then adjoint variant (Sect.~\ref{subsection:vanilla}). In step 1, the state equation is fully solved including forward communication; afterwards the adjoint equation is solved, including backward communication.}
\label{fig:vanilla}
\end{figure}
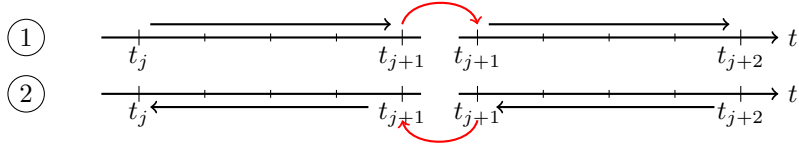

In this straightforward approach, PFASST is used first to solve the state equation~\eqref{eq:state1},~\eqref{eq:state2}. The state solution is recorded at the quadrature nodes, and then used in a subsequent PFASST run to solve the adjoint equation~\eqref{eq:adjoint1},~\eqref{eq:adjoint2}, see Fig.~\ref{fig:vanilla} for a sketch. This has already been used in the preliminary study~\cite{GoetschelMinion2018}.

Compared to sequential time stepping, the storage overhead is larger here, as the discretized-in-space state solution has to be stored on quadrature nodes, not only on the time steps. This is somewhat mitigated by the fact that the time intervals are distributed across several processors, typically giving access to more memory. Using PFASST gives some flexibility in reducing memory requirements. Besides compressed storage, e.g.,~\cite{GoetschelWeiser2015} for adaptive lossy compression of finite element solutions, options include storing the state on the coarse level only, and reducing the compression error by additional state sweeps without communication for the adjoint solve. A detailed analysis of storage requirements and storage reduction techniques will be reported elsewhere.

\subsubsection{Simultaneous approach}\label{subsection:simul}

Alternatively, if $R=N$, the adjoint can be solved simultaneously with the state equation (Fig.~\ref{fig:combined}), requiring communication of updated initial values for the state equation forward-in-time as well as updated terminal conditions for the adjoint backward-in-time. This cross-communication makes an efficient implementation with overlapping communication and computation difficult (see Fig.~\ref{fig:combined_pfasst}). As our numerical experiments show that this induces severe wait times, thus rendering the method inefficient, we do not consider this approach further.

\begin{figure}
\begin{center}
\begin{tikzpicture}
\node[shape=circle,draw,inner sep=2pt] at (-1,0) (char){1};
\draw[thick, black] (0,0) -- (4.25,0);
\draw[thick, black, ->] (4.75,0) -- (9,0) node[right] {$t$};
\draw (0.5,-3pt) -- (0.5,3pt) node[below = 1mm] {$t_{j}$};
\draw (1.375,-1.5pt) -- (1.375,1.5pt) ;
\draw (2.25,-1.5pt) -- (2.25,1.5pt) ;
\draw (3.125,-1.5pt) -- (3.125,1.5pt) ;
\draw (4.0,-3pt) -- (4.0,3pt) node[below = 1mm] {$t_{j+1}$};
\draw (5.0,-3pt) -- (5.0,3pt) node[below = 1mm] {$t_{j+1}$};
\draw[thick, black, ->] (0.65, 5pt) -- (3.85, 5pt);
\draw[thick, red, ->] (4.0, 5pt) to [out=75, in=105] (5.0, 5pt); 
\draw[thick, black, ->] (5.15, 5pt) -- (8.35, 5pt);
\draw (5.875,-1.5pt) -- (5.875,1.5pt) ;
\draw (6.75,-1.5pt) -- (6.75,1.5pt) ;
\draw (7.625,-1.5pt) -- (7.625,1.5pt) ;
\draw (8.5,-3pt) -- (8.5,3pt) node[below = 1mm] {$t_{j+2}$};
\draw[thick, black, <-] ($(0.65, 0)+(0,-5pt)$) -- ($(3.55,0)+(0,-5pt)$);
\draw[thick, red, <-] ($(4.0, 0)+(0,-10pt)$) to [out=-75, in=-105] ($(5.0,0)+(0,-10pt)$); 
\draw[thick, black, <-] ($(5.25, 0)+(0,-5pt)$) -- ($(8.15,0)+(0,-5pt)$);
\end{tikzpicture}
\end{center}
\caption{Simultaneous solve of state and adjoint (Sect.~\ref{subsection:simul}), requiring forward and backward communication between two adjacent timesteps.}
\label{fig:combined}
\end{figure}
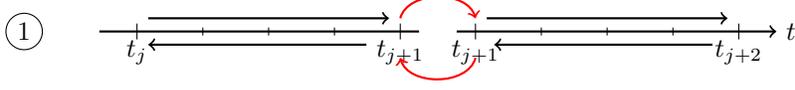
\begin{figure}
\begin{center}
\includegraphics[width=0.9\textwidth]{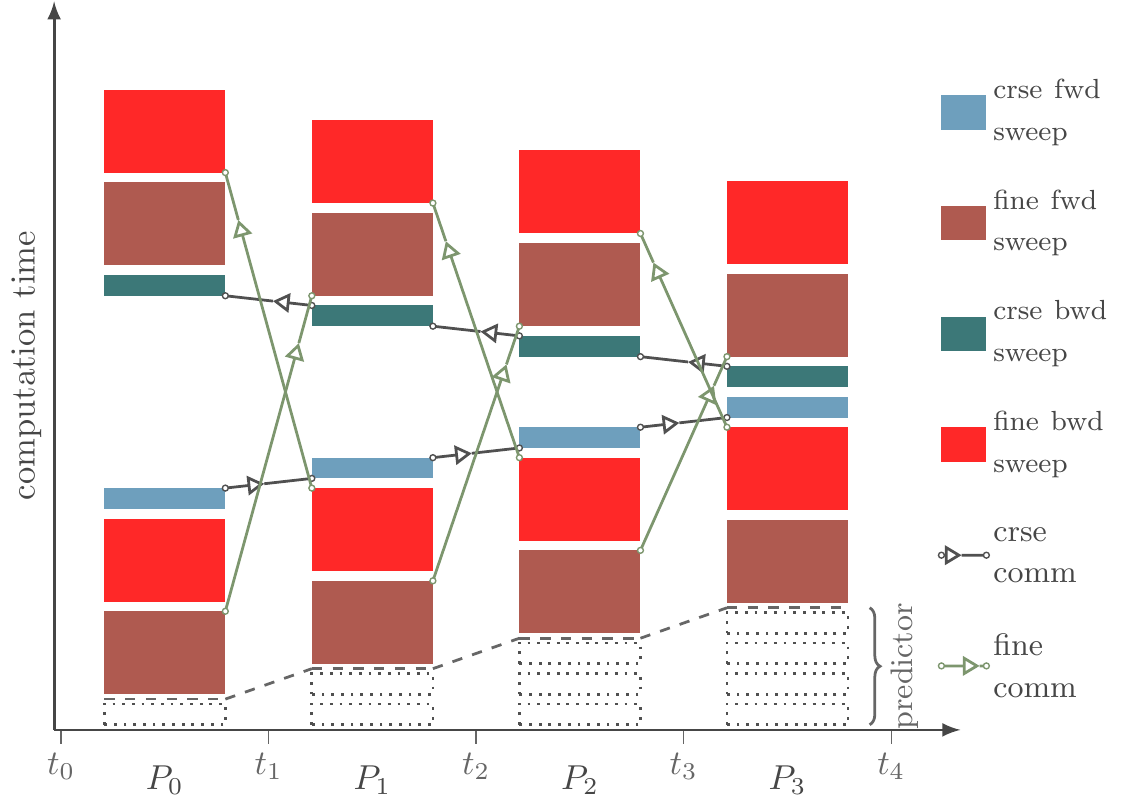}
\caption{Sweeps and communication pattern for simultaneously solving state and adjoint equations. The forward-backward-communication makes an efficient implementation without severe wait times difficult. Picture created using a modified version of pfasst-tikz (\protect\url{https://github.com/f-koehler/pfasst-tikz}).}
\label{fig:combined_pfasst}
\end{center}
\end{figure}

\subsubsection{Mixed approach}\label{subsection:mixed}

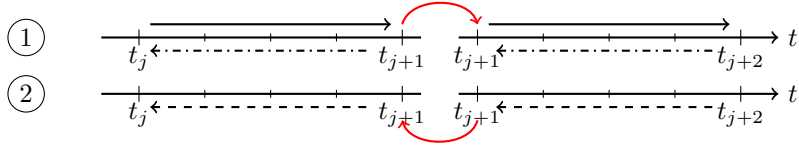
\begin{figure}
\begin{center}
\begin{tikzpicture}
\node[shape=circle,draw,inner sep=2pt] at (-1,0) (char){1};
\draw[thick, black] (0,0) -- (4.25,0);
\draw[thick, black, ->] (4.75,0) -- (9,0) node[right] {$t$};
\draw (0.5,-3pt) -- (0.5,3pt) node[below = 1mm] {$t_{j}$};
\draw (1.375,-1.5pt) -- (1.375,1.5pt) ;
\draw (2.25,-1.5pt) -- (2.25,1.5pt) ;
\draw (3.125,-1.5pt) -- (3.125,1.5pt) ;
\draw (4.0,-3pt) -- (4.0,3pt) node[below = 1mm] {$t_{j+1}$};
\draw (5.0,-3pt) -- (5.0,3pt) node[below = 1mm] {$t_{j+1}$};
\draw[thick, black, ->] (0.65, 5pt) -- (3.85, 5pt);
\draw[thick, red, ->] (4.0, 5pt) to [out=75, in=105] (5.0, 5pt); 
\draw[thick, black, ->] (5.15, 5pt) -- (8.35, 5pt);
\draw (5.875,-1.5pt) -- (5.875,1.5pt) ;
\draw (6.75,-1.5pt) -- (6.75,1.5pt) ;
\draw (7.625,-1.5pt) -- (7.625,1.5pt) ;
\draw[thick, black, dash dot, <-] ($(0.65, 0)+(0,-5pt)$) -- ($(3.55,0)+(0,-5pt)$);
\draw[thick, black, dash dot, <-] ($(5.25, 0)+(0,-5pt)$) -- ($(8.15,0)+(0,-5pt)$);
\draw (8.5,-3pt) -- (8.5,3pt) node[below = 1mm] {$t_{j+2}$};
\node[shape=circle,draw,inner sep=2pt] at (-1,-0.75) (char){2};
\draw[thick, black] (0,-0.75) -- (4.25,-0.75);
\draw[thick, black, ->] (4.75,-0.75) -- (9,-0.75) node[right] {$t$};
\draw ($(0.5,-0.75)+(0,-3pt)$) -- ($(0.5,-0.75)+(0,+3pt)$) node[below = 1mm] {$t_{j}$};
\draw ($(1.375,-0.75)+(0,-1.5pt)$) -- ($(1.375,-0.75)+(0,1.5pt)$);
\draw ($(2.25,-0.75)+(0,-1.5pt)$) -- ($(2.25,-0.75)+(0,1.5pt)$) ;
\draw ($(3.125,-0.75)+(0,-1.5pt)$) -- ($(3.125,-0.75)+(0,1.5pt)$);
\draw ($(4.0,-0.75)+(0,-3pt)$) -- ($(4.0,-0.75)+(0,3pt)$) node[below = 1mm] {$t_{j+1}$};
\draw ($(5.0,-0.75)+(0,-3pt)$) -- ($(5.0,-0.75)+(0,3pt)$) node[below = 1mm] {$t_{j+1}$};
\draw[thick, black, dashed, <-] ($(0.65, -0.75)+(0,-5pt)$) -- ($(3.55,-0.75)+(0,-5pt)$);
\draw[thick, red, <-] ($(4.0, -0.75)+(0,-10pt)$) to [out=-75, in=-105] ($(5.0, -0.75)+(0,-10pt)$); 
\draw[thick, black, dashed, <-] ($(5.25, -0.75)+(0,-5pt)$) -- ($(8.15,-0.75)+(0,-5pt)$);
\draw ($(5.875,-0.75)+(0,-1.5pt)$) -- ($(5.875,-0.75)+(0,1.5pt)$) ;
\draw  ($(6.75,-0.75)+(0,-1.5pt)$) -- ($(6.75,-0.75)+(0,1.5pt)$) ;
\draw  ($(7.625,-0.75)+(0,-1.5pt)$) -- ($(7.625,-0.75)+(0,1.5pt)$) ;
\draw ($(8.5,-0.75)+(0,-3pt)$) -- ($(8.5,-0.75)+(0,3pt)$) node[below = 1mm] {$t_{j+2}$};
\end{tikzpicture}
\end{center}
\caption{Mixed approach (Sect.~\ref{subsection:mixed}). In step 1, the full state equation is solved forward in time including forward communication. Simultaneously, a modified adjoint equation is solved without backward communication. In step 2, a correction for the adjoint is computed, propagating the correct terminal values.}
\label{fig:mixed}
\end{figure}

Making use of the linearity of the adjoint equation, the adjoint $p$ described in
equation~\eqref{eq:adjoint1},~\eqref{eq:adjoint3} can split into two terms $p = \tilde p + \delta$, where  $\tilde p$ satisfies the same equation
as $p$ but with homogeneous terminal conditions, and the defect $\delta$ which corrects the terminal conditions.
Specifically, for a given state solution $\bar y$, let
\begin{alignat}{3}
-\tilde p_t - \kappa \Delta\tilde p + f_y(\bar y)\tilde p &= -J^\Omega(\bar y, u) &\quad&\text{in}\ \Omega\times [0,T]\label{eq:ptilde} \\
\tilde p (T) &= 0 &&\text{in}\ \Omega \label{eq:ptildeT}
\end{alignat}
and
\begin{alignat}{3}
-\delta_t - \kappa \Delta\delta + f_y(\bar y)\delta &= 0 &\quad&\text{in}\ \Omega\times [0,T]\label{eq:delta} \\
\delta (T) &= -J^T(\bar y,u) &&\text{in}\ \Omega.\label{eq:deltaT}
\end{alignat}
To allow for an efficient time-parallel solution, we have to further modify the equations on the time intervals.
Denote by the superscript $j$, e.g.,~$\tilde{p}^j, \delta^j$, the respective solution on the $j$th time interval $[t_j, t_{j+1}]$.
With this notation, let $\tilde p^j$ solve
\begin{alignat}{3}
-\tilde p^j_t - \kappa \Delta\tilde p^j + f_y(\bar y)\tilde p^j &= -J^\Omega(\bar y, u) &\quad&\text{in}\ \Omega\times [t_j, t_{j+1}]\label{eq:adjoint_mod_interval1} \\
\tilde p^j (t_{j+1}) &= 0 &&\text{in}\ \Omega.\label{eq:adjoint_mod_interval2}
\end{alignat}
Note the homogeneous boundary condition \emph{on the interval}. For the defect $\delta^j$, the respective equation then is
\begin{equation}\label{eq:defect1}
-\delta^j_t - \kappa \Delta\delta^j + f_y(\bar y)\delta^j = 0 \quad\text{in}\ \Omega\times [t_j,t_{j+1}]
\end{equation}
with terminal conditions
\begin{equation}\label{eq:defect2}
\delta^j (t_{j+1}) = \begin{cases} \tilde p^{j+1}(t_{j+1}) + \delta^{j+1}(t_{j+1})  &j = 0,\dots, N-1\\
									-J^T(\bar y,u) &j=N.
					\end{cases}
\end{equation}
With this, equations~\eqref{eq:adjoint_mod_interval1}, \eqref{eq:adjoint_mod_interval2} can be solved in parallel together with the state equation on the interval without requiring communication backwards-in-time (first step in Fig.~\ref{fig:mixed}, the dash-dotted line denotes the modified adjoint equation). In the second step, the defect equation is solved (dashed line in Fig.~\ref{fig:mixed}). Note that on each time interval, the adjoint solvers only need to access the state solution computed on the same interval (even if each processor computes multiple time steps), thus requiring no additional communication of state values.

\begin{remark}
In case of no terminal cost in the objective function, the splitting of the adjoint is still required, since the homogeneous terminal condition only holds at time $t=T$. After time discretization, the terminal condition \emph{of the timestep}, i.e., $\tilde p^j (t_{j+1})$, is in general non-zero. To avoid backward-in-time communication these terminal conditions on the respective intervals are treated by the defect equation~\eqref{eq:defect1}, \eqref{eq:defect2}.
\end{remark}

For the solution of~\eqref{eq:defect1}, \eqref{eq:defect2} we note that formally, defining the operator $L = -\kappa \Delta + f_y(\bar y) ,$ the solution is given by
\begin{equation}\label{eq:exponential}
\delta^j(t) = \exp\bigl((t-t_{j+1})L\bigr)\delta^{j+1} (t_{j+1}) \quad t\in[t_j,t_{j+1}].
\end{equation}

For linear state equations we have $f_y = 0$, and equation~\eqref{eq:exponential} can be solved efficiently using suitable approximations of the matrix exponential, see the discussion in~\cite{GanderGuettel2013}. Thus, for the numerical experiments we restrict ourselves to the linear case. Besides having a certain relevance on their own, such equations occur, e.g., as subproblems in an SQP method, or in the evaluation of Hessian-times-vector products for Newton-CG methods, where an additional linearized state equation and a corresponding second adjoint equation need to be solved. For nonlinear state equations, the differential operator and/or coefficients of the adjoint are time-dependent (as they depend on the state solution). The extension of the mixed approach to this case (e.g., facilitating the Magnus expansion along the lines of~\cite{KrullMinion2018a}, where temporal parallelism of Magnus integrators is explored) is left for future work.

\begin{remark}
The mixed approach is to some extent related to the paraexp algorithm~\cite{GanderGuettel2013} for linear initial value problems. While the treatment of the source term on time intervals $[t_j, t_{j+1}]$ is the same, the paraexp method uses a near-optimal exponential integrator to compute in parallel solutions with the correct initial value on intervals $[t_j,T], j=0,\dots,N-1$, i.e., up to the final time $T$, and then sums up the individual solutions, requiring communication between several processors for computing the superposition at the required time points. In contrast, we solve the defect equation to propagate the correct terminal values on the time interval $[t_j, t_{j+1}], j=0,\dots,N-1$ only, such that only solution values at $t_j,\, j=N,\dots,1$ have to be communicated to the adjacent processor.
\end{remark}

\subsection{Warm starts} \label{sect:warm}

An important benefit of using SDC sweeps is the option to re-use information from the previous optimization iteration to initialize the subsequent PDE solves. For this, the solution of state and/or adjoint at the quadrature nodes is stored in one optimization iteration, and used as an initial guess instead of the predictor step during the following state/adjoint solve with an updated control. As the optimization is converging, $\norm{u_{k+1}-u_k}$ gets smaller, so that the previous solution of the state and adjoint equations are a suitable initialization, and the PFASST algorithm will take fewer iterations to converge.  
This incurs an overhead in storage that can be reduced by either using compression, or, as it is only used as initialization, values can be stored only on the 
coarsest level, and interpolated.

For the numerical examples below, to sweep on state or adjoint in optimization iteration $k+1$ we load the stored fine-level solution of the previous iteration $k$ to initialize at the fine quadrature nodes. Instead of invoking the usual PFASST predictor, we restrict the solution values and function evaluations to the coarser levels, and perform sweeps on the coarse level with communication as described in Sect.~\ref{sect:pfasst}.

\section{Numerical Results}
\label{section:results}

In this section, we test the developed methods on two optimal control problems. First, we consider optimal control for a linear heat equation in Sect.~\ref{subsection:heat} to compare the first state, then adjoint approach and the mixed approach for cold and warm starts. As a second example, an optimal control problem governed by the nonlinear Nagumo reaction-diffusion 
equation is considered (Sect.~\ref{subsection:nagumo}). There we compare different methods of treating the nonlinearity and again consider cold and warm starting of the time integration, using the first state, then adjoint approach. All examples are implemented using the LibPFASST library~\cite{libpfasst}, and are timed on Intel Xeon E5-4640 CPUs clocked at 2.4GHz.

Before discussing the results in detail, two remarks are in order. First, we report speedup and parallel efficiency with respect to the sequential versions of the described methods, i.e., MLSDC for state and adjoint solution. In particular, the sequential integration method has the same temporal order as the parallel version. We do not perform a thorough comparison of our methods with other commonly used sequential temporal integrators, as it is not clear what the ``best'' sequential method is for the respective problems. However, in Sect.~\ref{subsection:nagumo}, we briefly report results using IMEX-Euler as well as a $4$th-order additive Runge-Kutta method as sequential references. Second, scaling experiments are performed for a relatively low number of processors (20 and 32 in the two examples). Usually, parallel-in-time methods are combined with parallelization in the spatial domain, and multiply the speedup (and the required number of processors). As a rule of thumb, time-parallel methods come into play when spatial parallelism saturates, i.e., most of the available processors will be used for spatial parallelism. 


\subsection{Linear problem: Heat equation}\label{subsection:heat}

For the first numerical test, we consider the linear-quadratic optimal control problem to minimize $J(y,u)$ given by~\eqref{eq:objective} subject to
\begin{alignat*}{3}
y_t - \Delta y &= u &\quad&\text{in}\ \Omega \times [0,T]\\
y(0) &= y_0 &&\text{in}\ \Omega
\end{alignat*}
with periodic boundary conditions. The spatial domain $\Omega = [0,1]^3$, with initial conditions
\begin{equation*}
y_0(x) = \frac{1}{12\pi^2\lambda}(1-T) \prod_{i=1}^3 \sin(2\pi x_i)\\ \quad x\in \Omega,
\end{equation*}
and target solution
\begin{equation*}
\begin{split}
y_d(x,t) = \biggl[\bigl(12\pi^2 + \frac{1}{12\pi^2\lambda}\bigr)(t-T)-\bigl(1+\frac{1}{(12\pi^2)^2\lambda}\bigr)\biggr] \prod_{i=1}^3 \sin(2\pi x_i),\\ x\in \Omega, t\in[0,T].
\end{split}
\end{equation*}
This implies the exact solution
\begin{align*}
p^\star &= (T-t) \prod_{i=1}^3 \sin(2\pi x_i) \\
u^\star &= -\frac{1}{\lambda}p^\star \\
y^\star &= \left(-\frac{1}{(12\pi^2)^2\alpha}\ t - \frac{1}{12\pi^2\alpha}\ T + \frac{1}{12\pi^2\lambda}\right) \prod_{i=1}^3 \sin(2\pi x_i).
\end{align*}


For strong scaling results we solve the above problem, with control cost parameter $\lambda=0.05$ in the objective \eqref{eq:objective}. We use 
a three level PFASST scheme, and a pseudo-spectral discretization in space with $(16/32/64)^3$ degrees of freedom on the respective levels.  The implicit linear solves are done via the FFT. 
In time, we use $T=2$ and $20$ time steps with  $2/3/5$ Lobatto IIIA quadrature rules, yielding a temporal order $8$. To solve the optimization problem we apply gradient descent with the Armijo step size rule, and stop the optimization after $50$ iterations to compare the optimization progress, computed controls, and timings of the different algorithmic variants (first state, then adjoint vs.~mixed, cold vs.~warm start). PFASST iterations are stopped when the absolute or relative residual drops below $10^{-10}$. We note that this residual tolerance is too strict for the initial optimization iterations, where less exact gradients would be sufficient. Future research will be concerned with a thorough analysis of accuracy requirements, and the influence of inexactness on the convergence of the optimization methods along the lines of~\cite{GoetschelWeiser2015}. A brief comparison of different residual tolerances is shown in Table~\ref{tab:heat:tolerances}. We note that while naturally the wallclock time decreases, the influence of increased residual tolerances on parallel efficiency is small, as both sequential and parallel versions are accelerated similarly.

Results are obtained using the first state, then adjoint approach as well as the mixed approach. For each variant we compare the usual PFASST predictor (cold start) to initialize the MLSDC sweepers on each time step, and warm starts, using the previously computed solution on the fine level to initialize the sweeper. Tables~\ref{tab:heat:scaling_plain_cold}--\ref{tab:heat:scaling_mixed_warm} show speedup and efficiency for strong scaling using up to $20$ processors for the different approaches. In Tables~\ref{tab:heat:scaling_plain_cold} and~\ref{tab:heat:scaling_mixed_cold}, the usual PFASST predictor is used, spreading the given initial value on each parallel time step to the quadrature nodes. While speedup in Table~\ref{tab:heat:scaling_plain_cold} shows the gain in computation time due to using PFASST instead of MLSDC, the speedup in Table~\ref{tab:heat:scaling_mixed_cold} is computed with respect to the sequential first state, then adjoint approach, i.e., contains speedup due to parallelism as well as splitting the adjoint solve. Tables~\ref{tab:heat:scaling_plain_warm} and~\ref{tab:heat:scaling_mixed_warm} show results for warm starting the respective methods. Total speedup is computed with respect to the sequential first state, then adjoint base method. Speedup to cold denotes the speedup compared to the same method and same number of processors, but using the standard PFASST predictor (burn-in instead of warm start, see Sect.~\ref{sect:pfasst}). In this example, warm starts increase the speedup by 10--20\%. For the first state, then adjoint variant, the 10\% gain due to warm start in the 20 processor setting corresponds to a reduction in MLSDC sweeps of 14\% (state) and 34\% (adjoint). Thus, for problems where the sweeps are more expensive, we expect warm starting to be more effective (compare also Sect.~\ref{subsection:nagumo}).

\begin{table}
\begin{center}
\pgfplotstabletypeset[timingsStyleHeat] 
{
nproc	walltime	speedup	speedup_to_cold
1	7808.9	0	1
2	6063.5	1.3	1
5	3797.3	2.1	1
10	2677.8	2.9	1
20	1679.9	4.6	1
}
\caption{Strong scaling results for the heat example: first state, then adjoint approach with cold start. Speedup and parallel efficiency are compared to the sequential run with one processor.}
\label{tab:heat:scaling_plain_cold}
\end{center}
\end{table}

\begin{table}
\begin{center}
\pgfplotstabletypeset[timingsStyleHeat] 
{
nproc	walltime	speedup		speedup_to_cold
1	8888.8		0		1
2	6116.0		1.3		1
5	3518.1		2.2		1
10	2508.1		3.1		1
20	1606.0		4.9		1
}
\caption{Strong scaling results for the heat example: mixed approach with cold start. Speedup computed with respect to the first state, then adjoint approach with one CPU.}
\label{tab:heat:scaling_mixed_cold}
\end{center}
\end{table}

\begin{table}
\begin{center}
\pgfplotstabletypeset[warmstartTimingsStyleHeat] 
{
nproc	walltime	speedup		speedup_to_cold
1	6474.1		1.2		1.2
2	5087.1		1.5		1.2
5	3172.1		2.5		1.2
10	2034.9		3.8		1.3
20	1347.5  	5.8		1.2
}
\caption{Strong scaling results for the heat example: first state, then adjoint with warm start. Total speedup computed with respect to first state, then adjoint approach with one processor and cold start. Speedup to cold denotes the speedup only due to warm start, i.e., compared to the time for cold start with the same number of CPUs. Efficiency is computed using the total speedup.}
\label{tab:heat:scaling_plain_warm}
\end{center}
\end{table}

\begin{table}
\begin{center}
\pgfplotstabletypeset[warmstartTimingsStyleHeat] 
{
nproc	walltime	speedup		speedup_to_cold
1	7187.6		1.1		1.2
2	5537.4		1.4		1.1
5	3424.8		2.3		1.0
10	2276.3		3.4		1.1
20	1494.6		5.2		1.1
}
\caption{Strong scaling results for the heat example: mixed approach with warm start. Total speedup computed with respect to the first state, then adjoint approach with one CPU and cold start. Speedup to cold denotes the speedup only due to warm start, i.e., compared to the time for the mixed approach using cold start with the same number of CPUs. Efficiency is computed using the total speedup.}
\label{tab:heat:scaling_mixed_warm}
\end{center}
\end{table}

\begin{table}
\begin{center}
\begin{tabular}{lcccc} \toprule
                         & \multicolumn{4}{c}{speedup (20 CPUs)} \\  
  residual tolerance     & \multicolumn{2}{c}{first state, then adjoint} & \multicolumn{2}{c}{mixed} \\
                   & cold & warm & cold & warm                    \\ \midrule
  $10^{-10}$  &  $4.9$ & $5.8$ & $4.9$ & $5.2$ \\
  $10^{-8}$   &  $4.9$ & $5.9$ & $4.6$ & $5.3$ \\
  $10^{-6}$   &  $4.5$ & $5.5$ & $4.3$ & $4.9$ \\
  $10^{-4}$   &  $4.1$ & $5.0$ & $3.9$ & $4.4$ \\ \bottomrule
\end{tabular}
\caption{Influence of residual tolerance on parallel speedup for the first state then adjoint approach with cold start. Note that the sequential reference is computed using the same tolerances, otherwise total speedups are computed as in Tables~\ref{tab:heat:scaling_plain_cold}--\ref{tab:heat:scaling_mixed_warm}.}
\label{tab:heat:tolerances}
\end{center}
\end{table}

\begin{figure}
\begin{center}
  \includegraphics{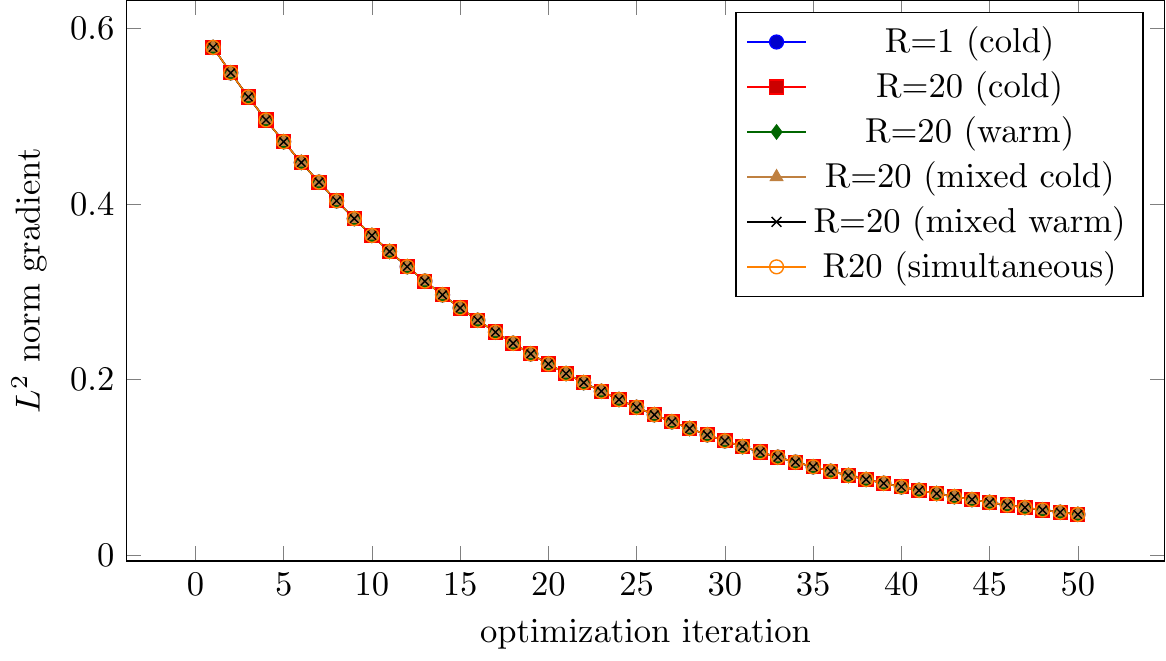}
  \caption{Heat example: Optimization progress is the same for the sequential first state, then adjoint approach (plain), its parallel variant (cold and warm start), and the parallel mixed approach (cold and warm start), as well as the simultaneous approach.}
  \label{fig:heat:optprogress}
  \end{center}
\end{figure}

\begin{figure}
\begin{center}
\includegraphics{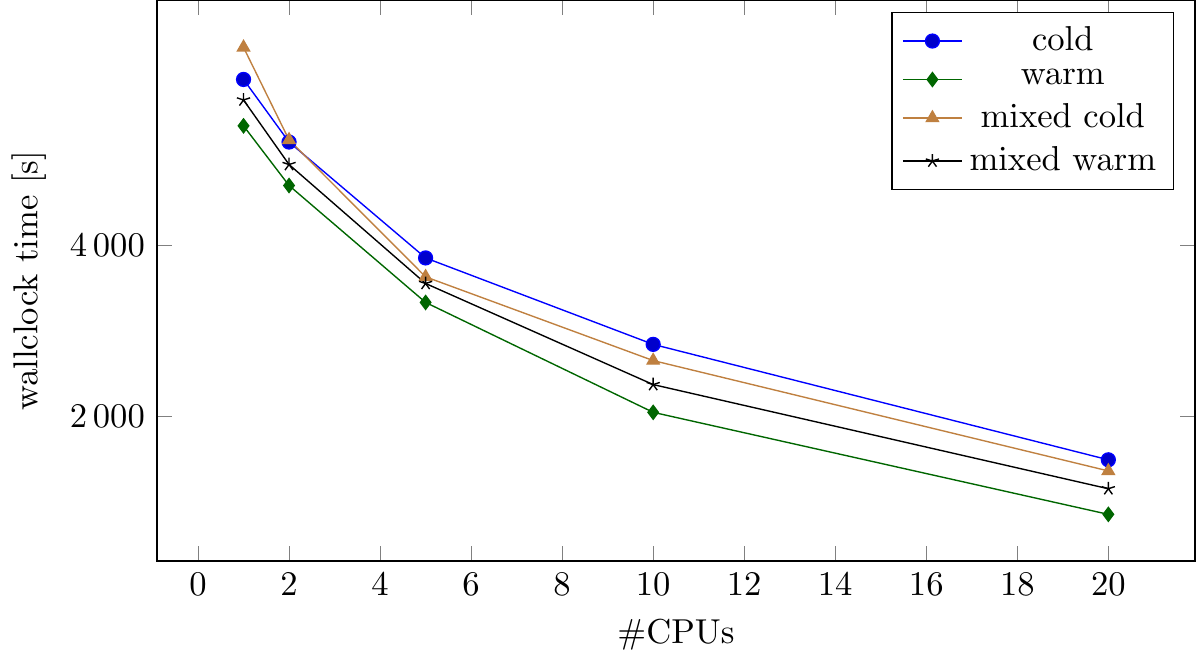}
\caption{Strong scaling for the heat example: comparison of wallclock times of mixed and first state, then adjoint approaches with cold and warm start.}
\label{fig:heat:compare}
\end{center}
\end{figure}

Fig.~\ref{fig:heat:optprogress} shows the progress of the optimization, which is similar for all approaches. The individual wallclock times are compared again in Fig.~\ref{fig:heat:compare}. For cold starting, the mixed approach performs better than the first state, then adjoint approach for most settings, but the advantage is rather small. 
Here, the mixed approach mainly saves computation time due to reduced communication, as both the propagation of the terminal condition in the mixed approach and the SDC sweeps in the plain first state, then adjoint approach are computed using the FFT and have a similar computational burden. Thus for other space discretizations, e.g., finite differences or finite elements, we expect the mixed approach to offer larger gains in computation times. 
Warm starting offers a moderate gain in speedup and efficiency for both variants (first state then adjoint, and mixed). As the gain is slightly higher in the first state, then adjoint approach, this variant provides the best performance in this test, achieving a parallel efficiency of $29\%$ on $20$ processors.

For completeness, the third variant, solving state and adjoint simultaneously (Sect.~\ref{subsection:simul}), leads to the same optimization progress (shown in Fig.~\ref{fig:heat:optprogress}), but gives a speedup of merely $1.8$ on $20$ processors, which is significantly worse than other variants. 

\subsection{Nonlinear problem: Nagumo equation}\label{subsection:nagumo}

Next we consider the following optimal control problem (see also~\cite{BuchholzEtAl2013, OPTPDE, OPTPDE2014}):
\begin{equation*}
 \min_{y,u} \half \int_0^T \int_\Omega (y-y_d)^2\ dx\ dt + \frac{\lambda}{2} \int_0^T \int_\Omega u^2\ dx\ dt,
\end{equation*}
i.e., as in the linear case, minimize~\eqref{eq:objective}, subject to
\begin{equation} \label{eq:nagumo}
 \begin{alignedat}{3}
    \frac{\partial}{\partial t}y(x,t) - \frac{\partial^2}{\partial x^2} y(x,t)  + (\frac{\gamma}{3} y^3(x,t) - y(x,t)) &= u(x,t) &\quad& \textrm{in}~\Omega \times (0,T)\\
    \frac{\partial}{\partial x}y(0,t) = \frac{\partial}{\partial x}y(L,t) &= 0 &&\textrm{in}~(0,T)\\
    y(x,0) &= y_0(x) && \textrm{in}~\Omega.
  \end{alignedat} 
\end{equation}
The parameter $\gamma$ steers the influence of the nonlinear reaction terms, and thus the stiffness due to the reaction term, and will be used below to compare IMEX and MISDC (multi-implicit SDC) formulations for solving the state equation. The spatial domain $\Omega = (0,L)$, with $L=20$, and the equation is solved up to final time $T = 5$. The initial condition is
\begin{equation*}
 y_0(x) = \begin{cases} 1.2\sqrt{3}, \quad &x \in [0,9] \\ 0, &\textrm{elsewhere}\end{cases}
\end{equation*}
and target solution
\begin{equation*}
 y_d(x,t) = \begin{cases} y_\text{nat}(x,t),\quad &t\in [0,2.5] \\ y_\text{nat}(x,2.5), &t \in (2.5,T],\end{cases}
\end{equation*}
where $y_\text{nat}$ denotes the solution to the PDE~\eqref{eq:nagumo} for $u \equiv 0$. Here, for $\lambda=0$, an exact optimal control is known,
\begin{equation*}
 u_\text{exact} = \begin{cases}
                   0, &t\leq 2.5\\
                   (\frac{\gamma}{3} y^3_\text{nat}(x,2.5) - y_\text{nat}(x,2.5)) - \frac{\partial^2}{\partial x^2} y_\text{nat}(x,2.5), \quad &t > 2.5,
                  \end{cases}
\end{equation*}
which will be used as a comparison for the computed optimal controls.
For this example, the adjoint equation is
\begin{equation}\label{eq:nagumo_adjoint}
 \begin{alignedat}{3}
    -\frac{\partial}{\partial t}p - \frac{\partial^2}{\partial x^2} p  + (\gamma y^2-1)p &= -(y-y_d) &\quad& \textrm{in}~\Omega \times (0,T)\\
    \frac{\partial}{\partial x}p(0,\cdot) = \frac{\partial}{\partial x}p(L,\cdot) &= 0 &&\textrm{in}~(0,T)\\
    p(\cdot,T) &= 0 && \textrm{in}~\Omega,
  \end{alignedat} 
\end{equation}
and the reduced gradient is given as $\nabla j(u) = \lambda u - p$.
For solving state and adjoint equation we do not consider the mixed method here, but always use the first state, then adjoint approach, since the adjoint  equation contains time-dependent terms in the differential operator (see the discussion in Sect.~\ref{subsection:mixed}).

\paragraph*{State equation and adjoint solve: MISDC vs.~IMEX}

Before coming to the actual optimization, let us briefly compare MISDC and IMEX-MLSDC for solving the state and adjoint equations~\eqref{eq:nagumo}, \eqref{eq:nagumo_adjoint} with zero control and varying $\gamma$. This demonstrates not only the flexibility of PFASST, but allows us to choose a suitable method in the scaling study for the optimization problem.

The IMEX and MISDC approach are explained by examining a single substep of an SDC sweep.
With $k$ denoting the SDC iteration,  $m$ the substep index, and
$\DD$ the discretization of the second derivative term, the 
fully implicit SDC version of an SDC substep \eqref{eq:sdcsweep}
for~\eqref{eq:nagumo} takes the form
\begin{equation}
  y^{[k+1]}_{m}=y_{j}  +  \qtil_{m,m}\Dt (\DD y^{[k+1]}_{m} -  y^{[k+1]}_{m} (\frac{\gamma}{3}(y^{[k+1]}_{m})^2-1))  +  S^{[k]}_m,
\end{equation}
where the term $S^{[k]}_m$ contains terms that either depend on the previous iteration $[k]$ or 
values  at  iteration $[k+1]$ already computed at substep $i<m$, 
including the control terms arising from the discretization of $u(x,t)$. 
The implicit equation couples nonlinear reaction and diffusion terms and hence requires
a global nonlinear solver in each substep.
For problems in which the reaction terms are non-stiff and can be treated explicitly,
the reaction terms at node $m$ do not appear in the implicit equation, giving the form
\begin{equation}
  y^{[k+1]}_{m}=y_{j}+\qtil_{m,m}\Dt (\DD y^{[k+1]}_{m})  +  S^{[k]}_m.
\end{equation}
Each IMEX substep now requires only the solution of a linear implicit equation, and thus
is computationally cheaper than the fully implicit approach, provided that the explicit treatment
of the reaction term does not impose an additional time step restriction.
When the reaction term is stiff, and hence it is advantageous to treat it implicitly,
a standard MISDC approach applies an operator splitting between diffusion and reaction in the
correction equation.  For example,
\begin{eqnarray}
  y^{*} & = & y_{j}+\qtil_{m,m}\Dt \DD y^* +  S^{*,[k]}_m, \label{eq:misdcI} \\
 y^{[k+1]}_{m} & = & y_{j}+\qtil_{m,m}\Dt (\DD y^{*} - y^{[k+1]}_{m} (\frac{\gamma}{3}(y^{[k+1]}_{m})^2-1))+ S^{[k]}_m. \label{eq:misdcII} 
\end{eqnarray}
We use Newton's method with damping and the natural monotonicity test according to~\cite{Deuflhard2006} to solve the nonlinear equation~\eqref{eq:misdcII} in each sweep. To reduce the computational effort, the MISDC approach is further modified so that
the nonlinear solve for reaction in~\eqref{eq:misdcII} is made linear by lagging terms in the splitting:
\begin{equation} \label{eq:lagging}
  y^{[k+1]}_{m}  = y_{j}+\qtil_{m,m}\Dt (\DD y^{*} - y^{[k+1]}_{m}(\frac{\gamma}{3}(y^{*})^2-1)) + S^{[k]}_m.
\end{equation}
This form creates an implicit solve with roughly the same cost as treating reaction explicitly but is
more stable.

\begin{figure}
 \centering
 \includegraphics[width=0.49\textwidth]{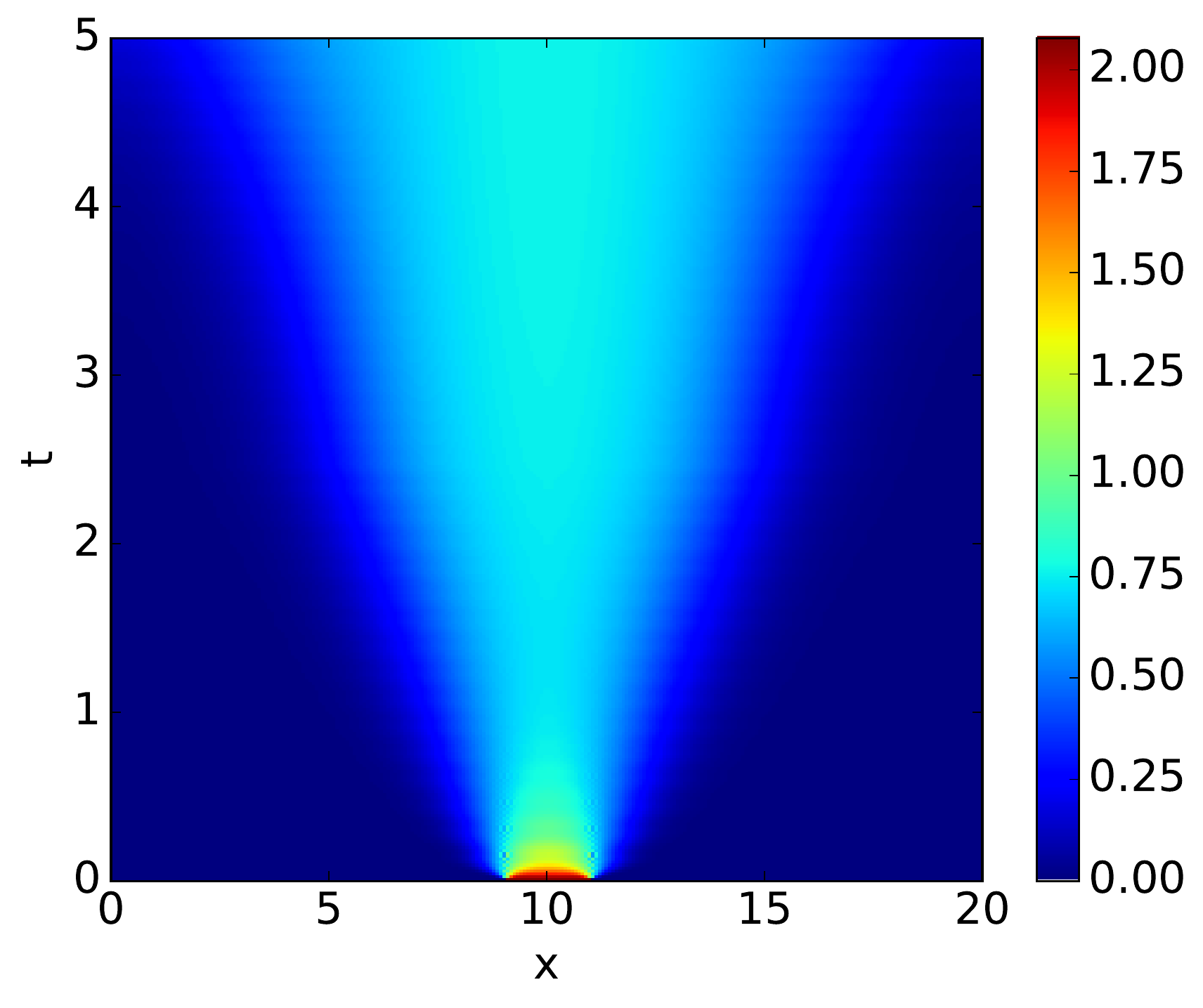}\ 
 \includegraphics[width=0.49\textwidth]{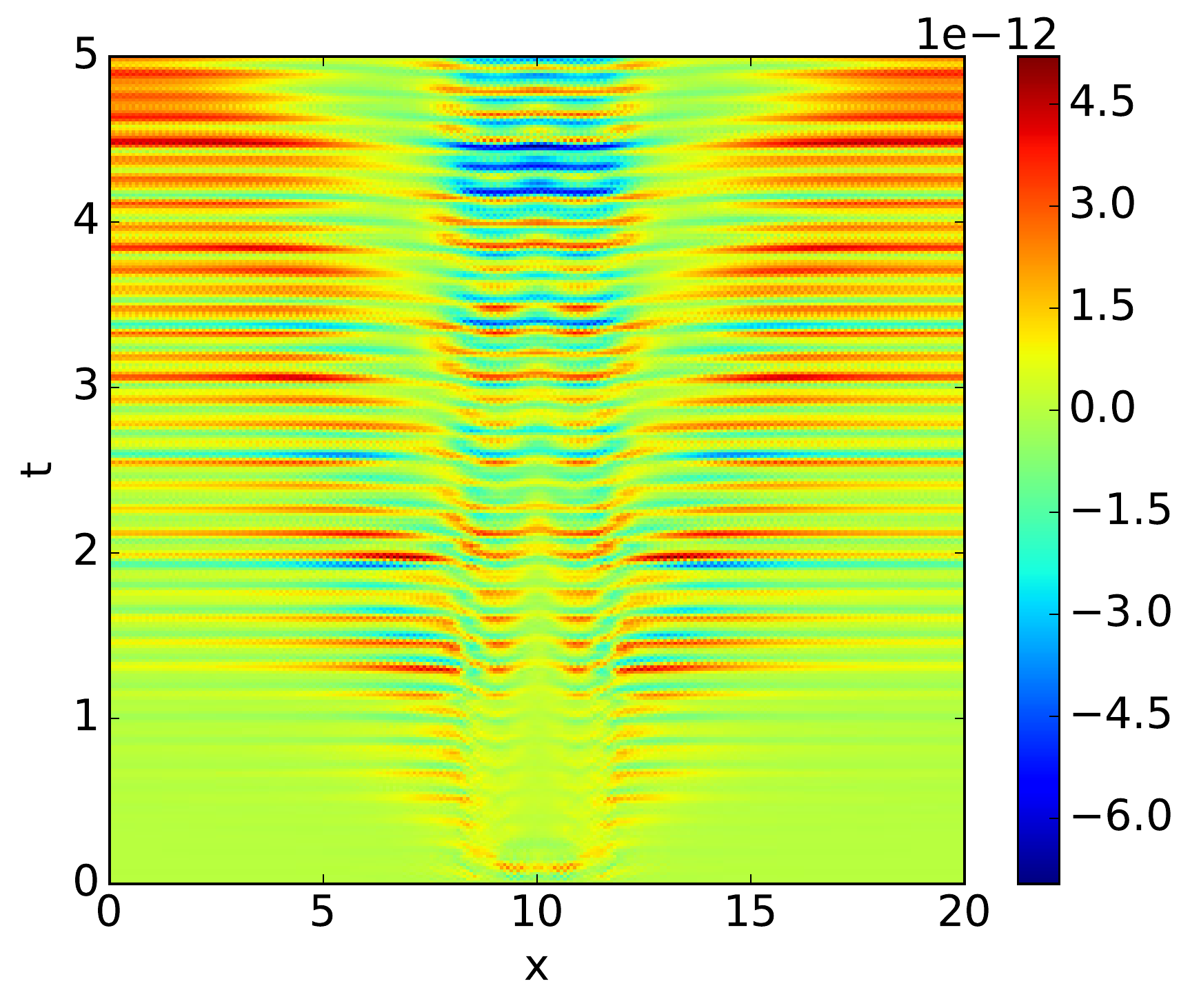}\ 
 \caption{Nagumo example: State solution for $\gamma=5$, $32$ parallel time steps using MISDC with lagging (left), and difference to the state solution obtained with nonlinear MISDC (right).}
 \label{fig:nagumo_misdc_sol}
\end{figure}

In all cases, a $3$-level PFASST scheme is applied to a  method of lines discretization. We use a pseudo-spectral discretization with 64/128/256 degrees of freedom in space, with spatial derivatives computed spectrally with the FFT. In time, $32$ uniform time steps are used, with 3/5/9 Lobatto IIIA collocation rules. Thus, all variants have the same order, and converge to the same collocation solution, if they converge. Sweeps are stopped when the relative or absolute residual drop below $10^{-11}$. Table~\ref{tab:nagumo:gamma} shows the required number of sweeps for convergence for several values of $\gamma$, time steps, and number of processors. While both MISDC variants stably converge for all values of $\gamma$ and time steps in sequential and parallel, IMEX requires finer time discretizations to converge, especially in parallel. In the sequential runs, the stability of the IMEX scheme requires a time step small enough so that the IMEX SDC iterations converge despite the explicit treatment of the reaction term. On the other hand, in PFASST, the sequential coarse grid IMEX sweeps are low-order and have a more restrictive stability constraint than the IMEX serial schemes.  When the solvers converge all variants require approximately the same number of sweeps on average, hence the IMEX scheme would be expected to need the least computational time. Fig.~\ref{fig:nagumo_misdc_sol} shows the solution for $\gamma=5$ and $32$ parallel time steps, computed with MISDC with lagging, and the difference to the solution obtained using MISDC with a nonlinear solver. While in most settings nonlinear and lagged versions need the same number of iterations, here the nonlinear variant requires on average $3$ iterations more than the lagged version, but the solutions agree to the convergence tolerance.

The conclusions of this comparison are two-fold.  First, the benefits of using an MISDC scheme versus an IMEX scheme will depend on the relative stiffness of the two terms and the relative cost of the nonlinear versus linear implicit equations.
Secondly,  at least for the problems considered here, the MISDC with lagging method performs as good or better than an IMEX splitting in terms of iteration count with only a small  computational overhead relative to the IMEX variant.

\begin{center}
\begin{table}
 \begin{tabular}{ccccccccc} \toprule
 \multirow{2}{*}{$\gamma$} & \multirow{2}{*}{$T/\Delta t$} & \multirow{2}{*}{\#CPUs} & \multicolumn{2}{c}{IMEX} & \multicolumn{2}{c}{MISDC~lagged} & \multicolumn{2}{c}{MISDC~nonlinear} \\
  & &   &  state & adjoint &  state & adjoint &  state & adjoint \\ \midrule
  \multirow{6}{*}{1} & \multirow{2}{*}{32}  & 1  & -   & -   & 14 & 9 & 14 & 9 \\
                     &                      & 32 & -   & -   & 53 & 23 & 55 & 23 \\ 
				     & \multirow{2}{*}{64}  & 1  & 7   & 6   & 7  & 6  & 7  & 6 \\
                     &                      & 32 & -   & -   & 31 & 19 & 32 & 19 \\ 
				     & \multirow{2}{*}{128} & 1  & 5   & 4   & 5  & 4  & 5  & 4 \\
                     &                      & 32 & 20  & 18  & 20 & 18 & 20 & 18 \\ \midrule
  \multirow{6}{*}{3} & \multirow{2}{*}{32}  & 1  & -   & -   & 14 & 9 & 14 & 9 \\
                     &                      & 32 & -   & -   & 53 & 24 & 54 & 24 \\ 
				     & \multirow{2}{*}{64}  & 1  & -   & -   & 7  & 6  & 7  & 6 \\
                     &                      & 32 & -   & -   & 31 & 19 & 32 & 19 \\ 
				     & \multirow{2}{*}{128} & 1  & 5   & 4   & 5  & 4  & 5  & 4 \\
                     &                      & 32 & 20  & 18  & 21 & 18 & 21 & 18 \\ \midrule
  \multirow{7}{*}{5} & \multirow{2}{*}{32}  & 1  & -   & -   & 14 & 9 & 14 & 9 \\
                     &                      & 32 & -   & -   & 53 & 24 & 56 & 24 \\ 
				     & \multirow{2}{*}{64}  & 1  & -   & -   & 7  & 6  & 7  & 6 \\
                     &                      & 32 & -   & -   & 31 & 19 & 32 & 19 \\ 
				     & \multirow{2}{*}{128} & 1  & 5   & 4   & 5  & 4  & 5  & 4 \\
                     &                      & 32 & -   & -   & 21 & 18 & 21 & 18  \\
					 & \multirow{2}{*}{256} & 1  & 4   & 4   & 4  & 4  & 4  & 4 \\
 					 &                      & 32 & 18  & 16  & 18 & 18 & 18 & 18 \\ \bottomrule					 			 
 \end{tabular}
 \caption{IMEX vs. MISDC for the Nagumo example: required number of PFASST iterations (average per time step, rounded). For the nonlinear solver, Newton's method is stopped when the discrete $\ell_2$-norm of the correction is smaller than $10^{-12}$.}
 \label{tab:nagumo:gamma}
\end{table}
\end{center}

\paragraph*{Optimization} To evaluate the performance of the parallel-in-time optimization method we use again $3$-level PFASST with an IMEX sweeper for solving state and adjoint equations up to a residual tolerance of $10^{-11}$. We set $\gamma=1$ in the state equation, and use 32/64/128 degrees of freedom in space as well as 32 time steps. In contrast to the experiments above, IMEX is converging in this setting due to the coarser space discretization, and has the smallest computation time. 
In the objective function the regularization parameter $\lambda=10^{-6}$ is used. Note that this leads to a difference between the computed optimal control and $u_\text{exact}$, as the latter is valid only for $\lambda=0$. The optimization is done using the nonlinear conjugate gradient method by Dai and Yuan~\cite{DaiYuan1999} in combination with linesearch to satisfy the the strong Wolfe conditions. As reported in~\cite{BuchholzEtAl2013}, the optimization progress with several ncg variants is slow, so we stop the optimization after 200 iterations. We again consider strong scaling with up to 32 CPUs. Due to the less severe nonlinearity in this test case and in view of the comparison above, only results using the IMEX variant are shown here. 

Tables~\ref{tab:nagumo:scaling_imex} and \ref{tab:nagumo:scaling_imex_warm} show timings, speedup, and parallel efficiency for IMEX with cold and warm starts. 
For cold starts, a speedup 
 of more than $5$ for $32$ processors is achieved, corresponding to a parallel efficiency of $16.5\%$.
\begin{table}
\begin{center}
\pgfplotstabletypeset[timingsStyle] 
{
nproc     walltime    speedup	speedup_to_cold	objective   final_rel_L2_err_ctrl   iter_state   iter_adjoint
1         169.8       0         1		5.7e-4      0.114                    293326      116565
2         109.0       1.6      	1		1.06e-3     0.127                    326573      125424
4         73.4        2.3       1		1.03e-3     0.125                    396014      153516
8         51.0        3.3       1 		4.0e-4      0.110                    529199      185789
16        37.7	      4.5       1		3.0e-4      0.105                    761860      247500
32        32.0        5.3       1		5.9e-4      0.115                   1086516      388938
}
\caption{Nagumo example: Strong scaling results for IMEX (first state, then adjoint; cold start). Speedup and parallel efficiency are compared to the sequential IMEX version with one CPU.}
\label{tab:nagumo:scaling_imex}
\end{center}
\end{table}

\begin{table}
\begin{center}
\pgfplotstabletypeset[warmstartTimingsStyle] 
{
nproc   walltime    speedup	speedup_to_cold  objective   final_rel_L2_err_ctrl   iter_state  iter_adjoint
1         132.7       1.3     	1.3          1.17e-3     0.130                   207642      117291
2         84.6        2.0		1.3	      	 1.16e-3     0.130                   229269      124787
4         58.0        2.9     	1.3       	 6.1e-4      0.115                   279609      148960
8         39.6        4.3     	1.3       	 7.5e-4      0.118                   355104      181228
16        30.4        5.6     	1.2  	   	 8.0e-4      0.120                   477929      242220
32        22.1        7.7    	1.4  	     5.3e-4      0.117                   660277		 338219
}
\caption{Nagumo example: Strong scaling results for IMEX (first state, then adjoint; warm start). Total speedup is compared to IMEX with one CPU and cold start, speedup to cold denotes the speedup only due to warm start, i.e., compared to the time for cold start with the same number of CPUs. Efficiency is computed using the total speedup.}
\label{tab:nagumo:scaling_imex_warm}
\end{center}
\end{table}

%

For the case with 32 CPUs, warm starting reduces the number of required state sweeps by $65\%$, while the reduction in adjoint sweeps is less pronounced. Overall, this translates to a decrease in computation time by nearly $40\%$, leading to a parallel efficiency of around $24\%$. The gain here is significantly larger than in the linear heat example (Sect.~\ref{subsection:heat}). This is true despite the tests being in one spatial dimension:  the PFASST method typically produces better parallel efficiencies for problems in higher dimensions since spatial coarsening then produces a greater relative reduction in computational cost on coarse levels.

Similar results (not reported here) are achieved using MISDC with lagging, albeit with overall higher computation times.

\begin{figure}
\begin{center}
  \includegraphics{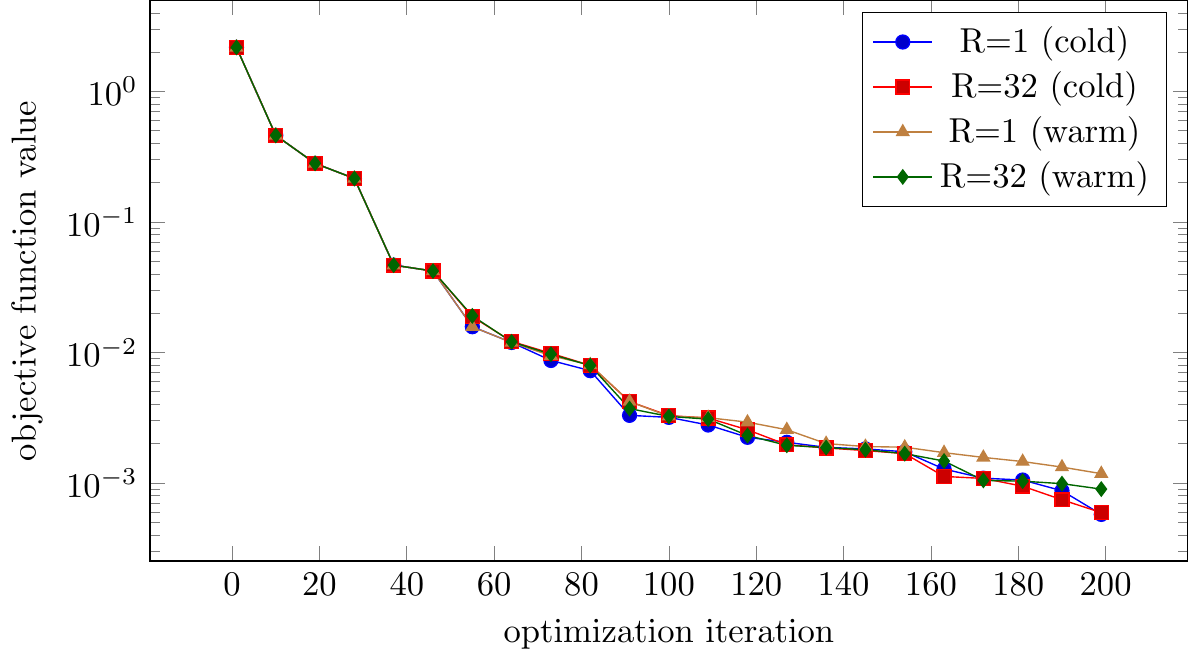}
  \caption{Nagumo example: Optimization progress for IMEX, using the first state, then adjoint approach, sequential and in parallel (32 CPUs), both with cold and warm starts.}
  \label{fig:nagumo:progress_imex}
  \end{center}
\end{figure}


Figure~\ref{fig:nagumo:progress_imex} 
shows the objective function value over the optimization process 
with 1 and 32 processors as well as using cold and warm starts.
We note that there are slight variances in the optimization progress and the computed controls, as also seen in the final error of the control in Table~\ref{tab:nagumo:scaling_imex}. The difference is that the sequential version does SDC sweeps on the later time steps with the correct initial value, while PFASST iteratively corrects the incorrect initial values. For warm starts, in addition the initial function evaluations as well as coarse grid corrections are inexact. Thus the iterative solution of the PDEs progresses differently, leading to slightly different final solutions and thus different gradients and step sizes, which accumulate during the optimization.

\paragraph{Comparison with classical time steppers}
The results reported above use MLSDC on one processor as the sequential reference for the computation of speedups. This allows the use of the same time discretization for all setups in the strong scaling study. However, due to the relatively high computational cost, MLSDC may not be the best sequential method. For comparison we consider a standard IMEX-Euler method as well as a fourth-order additive Runge-Kutta (ARK-4) method of Kennedy and Carpenter~\cite{KennedyCarpenter2003}. The IMEX-Euler method treats the reaction part explicitly, the diffusion part implicitly; ARK-4 combines an explicit Runge-Kutta method for the reaction part with an ESDIRK method for the diffusion. 
To obtain a reference solution for state and adjoint equations, we fix the spatial discretization to 128 degrees of freedom as before and use SDC with 9 Lobatto IIIA collocation nodes on a very fine temporal mesh. Time step sizes $\Delta t$ for the different methods are chosen such that the time discretization error compared to the reference solution is approximately the same for each case, which would result in similar progress of the full optimization process. All other parameters are chosen as before. 

One complication in this comparison is that computing the stage values in the Runge-Kutta method requires  evaluating the control $u$ and desired state $y_d$ at intermediate times (as does the MLSDC method require this data at the collocation nodes). For MLSDC we chose symmetric collocation nodes, such that forward-in-time and backward-in-time use the same nodes. Thus, $u, y_d$ have distinct values at the collocation nodes. However, in the ARK-4 method, the stage values are not symmetric, and interpolation or dense output would be required to compute these values. For simplicity, we therefore (artificially) set $u$ and $y_d$ constant on the time steps for this example.
Results are reported in table~\ref{tab:nagumo:seqref}.
\begin{table}
\begin{center}
\begin{tabular}{lrrr}\toprule
method      & $\Delta t$            & time [s] & rel. err. to $u_\text{exact}$ \\\midrule
IMEX-Euler  & $10^{-3}$             & $102.5$ & $1.2\cdot 10^{-1}$\\
ARK-4       & $2.5\cdot 10^{-3}$    & $183.0$ & $1.4\cdot 10^{-1}$\\
MLSDC cold  & $1.5625\cdot 10^{-1}$ & $169.8$ & $1.1\cdot 10^{-1}$\\
MLSDC warm  & $1.5625\cdot 10^{-1}$ & $132.7$ & $1.3\cdot 10^{-1}$\\\bottomrule
\end{tabular}
\caption{Runtimes of different sequential methods yielding similar optimization results.}
 \label{tab:nagumo:seqref}
\end{center}
\end{table}
Despite the small time step size, IMEX-Euler is the fastest of the sequential methods. However, using PFASST with warm starting, speedup is achieved already by adding the second processor (speedup $1.2$). Compared to IMEX-Euler, on $32$ processors, PFASST with cold starting achieved a speedup/parallel efficiency of $3.2$/$10\%$, with warm starts this increases to $4.6$/$14\%$.  As noted before, the benefit of MLSDC compared to SDC is greater for problems in more than one-dimension, hence these results should not be considered as a general comparison.

\section{Conclusions and Outlook} \label{section:outlook}

In this paper, we introduce an efficient fully time-parallel strategy for gradient-based optimization with parabolic PDEs.
In the most critical component, the computation of the reduced gradient, we discuss and compare three competing
approaches for applying the PFASST algorithm to the solution of the state and adjoint equations. 
While simultaneously solving state and adjoint induces severe communication and wait times in the implementation, both the first state, then adjoint approach and the mixed approach yield comparatively good speedups and parallel efficiency of up to $29\%$. The warm starting capability of MLSDC/PFASST is an important feature, and makes the method competitive with standard time stepping methods even in a sequential run. It is also noteworthy that speedup is obtained already with two processors. 

Although the presented results are promising, there are several additional avenues to further increase the parallel speedup or efficiency of the approach. For example, the test cases studied here were set in simple periodic domains using the FFT for implicit solves.  In many cases, implicit solves are done iteratively and hence solver tolerances can be adjusted dynamically to do less work when less accuracy is needed \cite{SpeckEtAl2016}.    The same is true of the tolerances used to terminate PFASST iterations.  
Likewise, the spatial or temporal order of the PDE solver could be changed dynamically during the optimization process.
Hence a dynamic, adaptive strategy for controlling these tolerances and accuracy in the optimization loop is a promising direction for future study.
Similarly, analyzing the impact of inexact storage of solution values for warm starting will be necessary for larger scale applications.  Finding suitable extensions of the mixed approach to nonlinear state equations is also of interest. For this approach, it might also be beneficial to solve state and modified adjoint equation in the first step \emph{concurrently}, making use of additional processors during SDC sweeps.

\subsection*{Acknowledgment}
  S.G. gratefully acknowledges partial funding by the Deutsche Forschungsgemeinschaft (DFG), Pro\-ject WE 2937/6-1.
  The work of M.M.  was supported by the Applied Mathematics Program of the
DOE Office of Advanced Scientific Computing Research under the U.S. Department
of Energy under contract DE-AC02-05CH11231, and the Alexander von Humboldt Foundation.


\end{document}